\pdfoutput=1
\RequirePackage{silence}
\documentclass[a4paper,svgnames]{amsart}
\usepackage[hmarginratio={1:1},vmarginratio={1:1},lmargin=60.0pt,tmargin=60.0pt]{geometry}

\usepackage{booktabs}
\allowdisplaybreaks

\usepackage[T1]{fontenc}
\usepackage{latexsym,exscale,amsfonts,amssymb,mathtools}
\usepackage{amsmath,amsthm,amsfonts,amssymb,amscd,textcomp,bbm}
\usepackage{mathrsfs,stackrel}
\usepackage{etoolbox}
\usepackage{tikz,tikz-cd}
\usetikzlibrary{matrix,arrows.meta,decorations.markings,shapes.multipart,decorations.pathreplacing,decorations.pathmorphing}
\usepackage{enumitem}
\setlist[enumerate]{itemsep=0.15cm,label=\emph{\upshape(\alph*)}}
\setlist[enumerate,2]{itemsep=0.15cm,label=\emph{\upshape(\roman*)}}
\usepackage{aliascnt}
\usepackage{ytableau}
\usepackage{xparse}
\usepackage{cite}

\usepackage{xcolor,colortbl}

\definecolor{mygray}{gray}{0.6}
\definecolor{mygraydark}{gray}{0.4}
\definecolor{mygraylight}{gray}{0.85}
\definecolor{spinach}{RGB}{46,139,87}
\definecolor{tomato}{RGB}{255,99,71}
\definecolor{orchid}{RGB}{143,40,194}
\definecolor{neon}{RGB}{77,77,255}
\definecolor{pumpkin}{RGB}{224,180,80}
\definecolor{citron}{RGB}{190,180,90}

\definecolor{lava}{RGB}{207,16,32}
\definecolor{cream}{RGB}{255,253,208}
\definecolor{verdigris}{RGB}{67,179,174}
\definecolor{Black}{RGB}{0,0,0}
\definecolor{mydarkblue}{RGB}{10,10,170}
\definecolor{darkspinach}{RGB}{20,70,20}
\definecolor{darktomato}{RGB}{155,40,30}
\definecolor{darkorchid}{RGB}{50,10,100}
\definecolor{darklava}{RGB}{150,8,16}

\usepackage{todonotes}

\newcommand{\BG}{{\mathbb G}}
\newcommand{\cD}{{\mathcal D}}
\newcommand{\cVect}[1][{\bf k}]{{\mathcal V}\mathrm{ect}_{#1}}
\newcommand{\sVect}[1][{\bf k}]{{\mathcal SV}\mathrm{ect}_{#1}}

\newcommand{\bK}{{\bf K}}
\newcommand{\bk}{{\bf k}}
\newcommand{\BZ}{{\mathbb Z}}
\newcommand{\BN}{{\mathbb N}}
\newcommand{\BQ}{{\mathbb Q}}
\newcommand{\BC}{{\mathbb C}}
\newcommand{\BR}{{\mathbb R}}
\newcommand{\BF}{{\mathbb F}}

\newcommand{\Ver}{\mathrm{Ver}}
\newcommand{\Rep}{\mathcal{R}\mathrm{ep}}
\newcommand{\coRep}{\mathrm{co}\mathcal{R}\mathrm{ep}}

\newcommand{\End}{\mathrm{End}}
\newcommand{\Aut}{\mathrm{Aut}}

\newcommand{\mN}{\mathbb{N}}
\newcommand{\mZ}{\mathbb{Z}}
\newcommand{\mC}{\mathbb{C}}
\newcommand{\ind}{\mathrm{ind}}
\newcommand{\res}{\mathrm{res}}
\newcommand{\mS}{\mathbb{S}}
\newcommand{\fg}{\mathfrak{g}}
\newcommand{\Lie}{\mathrm{Lie}\,}
\newcommand{\Dist}{\mathrm{Dist}\,}
\newcommand{\oa}{\bar{0}}
\newcommand{\ob}{\bar{1}}
\newcommand{\cO}{\mathcal{O}}

\newcommand{\placeholder}{{}_{-}}
\newcommand{\mystrut}{\rule[-0.2\baselineskip]{0pt}{0.9\baselineskip}}

\usepackage{tikz}
\tikzset{anchorbase/.style={baseline={([yshift=-0.5ex]current bounding box.center)}},}

\def\NewTheorem#1{%
\newaliascnt{#1}{equation}%
\newtheorem{#1}[#1]{#1}%
\aliascntresetthe{#1}%
\expandafter\def\csname #1autorefname\endcsname{#1}%
}
\def\equationautorefname~#1\null{(#1)\null}

\numberwithin{equation}{subsection}

\NewTheorem{Proposition}
\NewTheorem{Theorem}
\NewTheorem{Corollary}
\AtEndEnvironment{Corollary}{\null\hfill$\square$}%
\NewTheorem{Lemma}
\theoremstyle{definition}
\NewTheorem{Definition}
\NewTheorem{Notation}
\NewTheorem{Question}
\NewTheorem{Example}
\AtEndEnvironment{Example}{\null\hfill$\Diamond$}%
\NewTheorem{Examples}
\AtEndEnvironment{Examples}{\vskip-10mm\null\hfill$\Diamond$}%

\theoremstyle{remark}
\NewTheorem{Remark}
\NewTheorem{Assumption}

\usepackage[hypertexnames=false]{hyperref}
\usepackage{bookmark}
\hypersetup{
pdftoolbar=true,
pdfmenubar=true,
pdffitwindow=false,
pdfstartview={FitH},
pdftitle={Growth rates of the number of indecomposable summands in tensor powers},
pdfauthor={Kevin Coulembier, Victor Ostrik and Daniel Tubbenhauer},
pdfsubject={},
pdfcreator={Kevin Coulembier, Victor Ostrik and Daniel Tubbenhauer},
pdfproducer={Kevin Coulembier, Victor Ostrik and Daniel Tubbenhauer},
pdfkeywords={},
pdfnewwindow=true,
colorlinks=true,
linkcolor=mydarkblue,
citecolor=teal,
filecolor=magenta,
urlcolor=orchid,
linkbordercolor=lava,
citebordercolor=teal,
urlbordercolor=orchid,
linktocpage=true
}

\newcommand{\nnfootnote}[1]{%
\begin{NoHyper}
\renewcommand\thefootnote{}\footnote{#1}%
\addtocounter{footnote}{-1}%
\end{NoHyper}
}

\def\makeautorefname#1#2{\csdef{#1autorefname}{#2}}
\makeautorefname{section}{Section}%
\makeautorefname{subsection}{Section}%
\makeautorefname{subsubsection}{Section}%

\begin{document}
\title[Growth rates of the number of indecomposable summands in tensor powers]{Growth rates of the number of indecomposable summands in tensor powers}
\author[K. Coulembier, V. Ostrik and D. Tubbenhauer]{Kevin Coulembier, Victor Ostrik and Daniel Tubbenhauer}

\address{K.C.: The University of Sydney, School of Mathematics and Statistics F07, Office Carslaw 717, NSW 2006, Australia}
\email{kevin.coulembier@sydney.edu.au}

\address{V.O.: University of Oregon, Department of Mathematics, Eugene, OR 97403, USA}
\email{vostrik@math.uoregon.edu}

\address{D.T.: The University of Sydney, School of Mathematics and Statistics F07, Office Carslaw 827, NSW 2006, Australia, \href{http://www.dtubbenhauer.com}{www.dtubbenhauer.com}, 
\href{https://orcid.org/0000-0001-7265-5047}{ORCID 0000-0001-7265-5047}}
\email{daniel.tubbenhauer@sydney.edu.au}

\begin{abstract}
In this paper we study the asymptotic behavior of the number of summands in tensor products of finite dimensional representations of affine (semi)group (super)schemes and related objects.
\end{abstract}

\nnfootnote{\textit{Mathematics Subject Classification 2020.} Primary: 
17B10, 18M05; Secondary: 16T05, 17B37, 20C25.}
\nnfootnote{\textit{Keywords.} Tensor products, asymptotic behavior, affine group schemes, affine semigroup schemes, semigroups, supergroups, Hopf algebras, (symmetric) monoidal categories.}

\addtocontents{toc}{\protect\setcounter{tocdepth}{1}}

\maketitle

\tableofcontents

\ytableausetup{centertableaux,mathmode,boxsize=0.5cm}

%%%%%%%%%%%%%%%%%%%%%%%%%%%%%%%%%%%%%%%%%%%%%%

\section{Introduction and main results}\label{s:main}

%%%%%%%%%%%%%%%%%%%%%%%%%%%%%%%%%%%%%%%%%%%%%%

A central, yet hard, problem in representation theory is the decomposition of tensor products of representations into indecomposable summands. Computations of these decomposition numbers are often major unsolved problems in representation theory.
In this paper we take a different perspective and we are interested in asymptotic properties of the number of indecomposables in tensor products 
of representation rather than explicit decompositions. In contrast to 
the question of explicit decompositions, we obtain results in extensive generality.

To get started,  let $\Gamma$ be a finite group with a finite dimensional representation $V$ over some field $\bk$.

\begin{Definition}\label{d:main}
We define
\begin{gather*}
b_{n}^{\Gamma,V}:=\#\text{indecomposable summands in $V^{\otimes n}$ counted with multiplicities}
\end{gather*}
(sometimes simply denoted by $b_n$ when $\Gamma$ and $V$ are clear from the context), where `indecomposable' means as $\Gamma$-representation. Let further
\begin{gather*}
\beta^{\Gamma,V}:=\lim_{n\to\infty}\sqrt[n]{b_{n}^{\Gamma,V}}.
\end{gather*}
\end{Definition}

Note that $b_{n}^{\Gamma,V}b_{m}^{\Gamma,V}\leq b_{n+m}^{\Gamma,V}$, so that $\beta^{\Gamma,V}$ is well-defined by (a version of) Fekete's Subadditive Lemma.

There is no chance 
to compute $b_{n}^{\Gamma,V}$ explicitly in this generality, but it turns out that the `limit' $\beta^{\Gamma,V}$ can be understood.
As a first step we note that we have the following lemma whose the proof 
is immediate.

\begin{Lemma}\label{l:main}
We have $b_{n}^{\Gamma,V}\leq(\dim V)^{n}$, and consequently
\begin{gather*}
\pushQED{\qed} 
\beta^{\Gamma,V}\leq\dim V.\qedhere
\popQED
\end{gather*}
\end{Lemma}

A classical result of Bryant--Kov{\'a}cs \cite[Theorem 1]{BrKo-tensor-group} shows that there exists some $n$ for which $V^{\otimes n}$ contains a projective direct summand (projective over $\Gamma$ if $V$ is a faithful and projective over the appropriate factor group of $\Gamma$ otherwise). As observed in \cite[Theorem 6.3]{BeSy-non-projective-part}, a consequence of this is that the bound for $\beta^{\Gamma,V}$ in \autoref{l:main} is actually an equality, that is:

\begin{Theorem}\label{t:main}
We have 
\begin{gather*}
\pushQED{\qed} 
\beta^{\Gamma,V}=\dim V.\qedhere
\popQED
\end{gather*}
\end{Theorem}

It is a natural question to which extent this result generalizes. Concretely, one can let $\Gamma$ be an infinite group, an infinite semigroup, an affine group scheme, super versions, a finite dimensional Hopf algebra or any other algebraic structure for which we have a notion of tensor products of representations, and $V$ a finite dimensional $\Gamma$-representation. \autoref{d:main} and \autoref{l:main} work verbatim, including that $\beta^{\Gamma,V}$ is well-defined, although clearly the above proof of \autoref{t:main} does not extend. 

We show that the theorem remains valid in great generality, but has limitations:

\begin{Theorem}\label{t2:main}
\leavevmode
\begin{enumerate}

\item \autoref{t:main} holds true for $\Gamma$ any affine semigroup superscheme as defined in \autoref{d:schemes}. \autoref{t:main} holds 
also true for the quantum groups $U_{q}(\mathfrak{sl}_{M})$ and $U_{q}(\mathfrak{gl}_{M})$ {and the vector representation}, where we allow any $q\in\bk^{\ast}$ for $\bk^{\ast}=\bk\setminus\{0,-1\}$ if $\mathrm{char}(\bk)\neq 2$ and 
$\bk^{\ast}=\bk\setminus\{0\}$ otherwise.

\item \autoref{t:main} does not hold true in general for {(comodules over)} Hopf algebras.

\end{enumerate}
\end{Theorem}

\begin{Remark}\label{r2:main}
\leavevmode
\begin{enumerate}

\item Note that \autoref{t2:main}.(a) includes the cases where 
$\Gamma$ is any, possibly infinite, abstract semigroup (for instance a monoid or group -- as observed above, \autoref{t2:main}.(a) is classical for finite groups, but already for
infinite abstract groups it appears to be new). Indeed, using Tannakian reconstruction we can associate an affine semigroup scheme to an abstract semigroup which has equivalent representation theory. The result can of course also be obtained directly without this observation and we elaborate 
on special cases like these in the main body of the paper.

\item We have not sought completeness in \autoref{t2:main}(a). We can, for example, include algebraic objects like transitive group{\em oids} in schemes, see \autoref{CorGS}. We omit such cases here for clarity and because they are included in \autoref{t3:main}(a) below.

\item The proof of \autoref{t2:main}.(a), excluding the quantum group case, can be reduced to the specific cases of the general linear group $GL_{M}$ and the general linear supergroup $GL_{M|N}$.

\item While proving \autoref{t2:main}.(a) we will also give asymptotic formulas in some cases, which improves \autoref{t2:main}.(a).

\end{enumerate}
\end{Remark}

Let us also comment on some variations of \autoref{t2:main}.(a) 
which have appeared in the literature.

\begin{Remark}\label{r2:main2}
\leavevmode
\begin{enumerate}

\item There are known variants of $\beta^{\Gamma,V}$ for which the analog of \autoref{t2:main}.(a) is not valid. For instance if, in characteristic $p>0$, one only counts direct summands of dimension not divisible by $p$ (categorical dimension zero), it is shown in \cite[Theorem 8.1]{CEO} for affine group superschemes that the corresponding limit yields a ring homomorphism from the Grothendieck ring which takes values in the ring of integers of a particular cyclotomic extension of $\BQ$. For example, for $p=5$, $\Gamma=\BZ/5\BZ$ and $V$ its indecomposable of dimension three the limit is the golden ratio. In the case of a finite group where one only counts non-projective summands one obtains the variant studied for instance in \cite{BeSy-non-projective-part}. (Note that these results do not cover the monoid and semigroup case and it would be interesting to know whether there 
are similar statements for monoids and semigroups.) Our result that $\beta^{\Gamma,V}=\dim V$
for affine semigroup superschemes might be useful to compute the 
various variants of $\beta^{\Gamma,V}$ in the literature.

\item Let $\bk=\overline{\bk}$ and consider $TL_{n}=\End_{SL_{2}}\big((\bk^{2})^{\otimes n}\big)$ which is known as the Temperley--Lieb algebra. By Schur--Weyl duality, the dimensions of simple $TL_{n}$-representations correspond to the 
decomposition multiplicities of indecomposables in $(\bk^{2})^{\otimes n}$, and thus, \autoref{t2:main}.(a) and its proof given below imply that some of the simple $TL_{n}$-representations are very large for $n\gg 1$.
In the spirit of this example, the paper \cite{KhSiTu-monoidal-cryptography} proposes to study large simple representations in the setting of monoids that arise in a `Schur--Weyl dual way' from \autoref{t2:main}.(a) such as the Temperley--Lieb monoid, the Brauer monoid etc. We hope that \autoref{t2:main}.(a) can be used to 
generalize \cite{KhSiTu-monoidal-cryptography} beyond the case of the monoids discussed in that work.

\end{enumerate}
\end{Remark}

Note that \autoref{t:main} is a weaker statement than 
an asymptotic formula for the growth of $b_{n}^{\Gamma,V}$.
We show this in the following example illustrating \autoref{t:main}, which is in some sense the key example.

\begin{Example}\label{e:main}
Take $\bk=\BC$, $\Gamma_{2}=SL_{2}$ and $V_{2}=\BC^{2}$, its vector representation. Then the first numbers $b_{n}^{\Gamma_{2},V_{2}}$ are:
\begin{gather*}
\{1,1,2,3,6,10,20,35,70,126,252\},\quad b_{n}^{\Gamma_{2},V_{2}}\text{ for }n=0,\dots,10.
\end{gather*}
Here, as throughout, we view $b_n=b_{n}^{\Gamma_{2},V_{2}}$ 
as a function in $n$.

A Mathematica loglog plot of $\sqrt[n]{b_{n}}$ (y-axis) for $n\in\{1,\dots,1000\}$ (x-axis) gives
\begin{gather*}
\begin{tikzpicture}[anchorbase]
\node at (0,0) {\includegraphics[height=5cm]{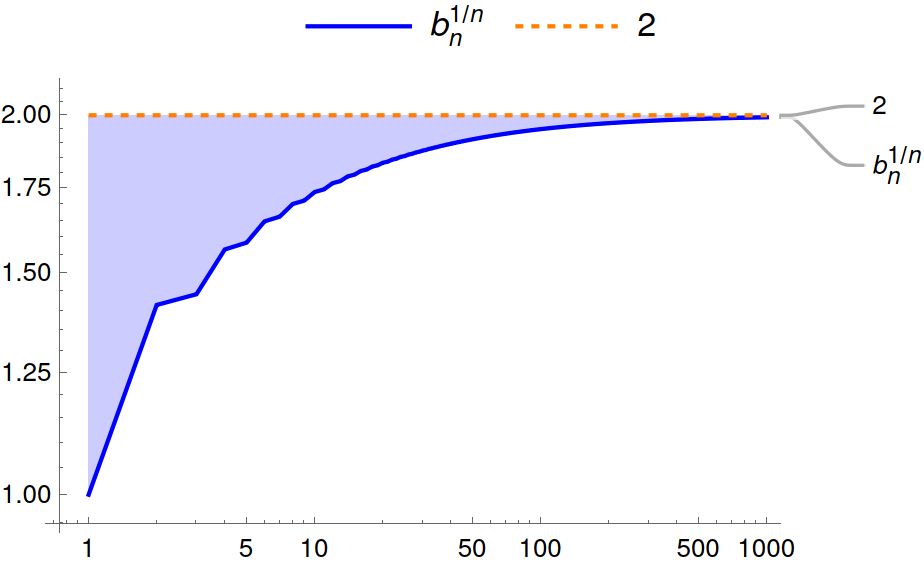}};
\end{tikzpicture}
,
\end{gather*}
and indeed the limit is two, as predicted by \autoref{t:main}. Precisely, 
$\sqrt[1000]{b_{1000}}\approx 1.99265$.

However, the asymptotic growth rate of $b_{n}$ is different than $2^{n}$.  As we will see in \autoref{e:clt},
we get 
\begin{gather*}
b_{n}^{\Gamma_{2},V_{2}}\sim\sqrt{2/\pi}\cdot 2^{n}/\sqrt{n}.
\end{gather*}
(Here and throughout, we use $f\sim g$ for $f$ is equal to $g$ asymptotically, meaning the ratio of $f$ and $g$ converges to one.)
We have $\sqrt{2/\pi}\approx 0.798$ and Mathematica's log plot gives:
\begin{gather*}
\begin{tikzpicture}[anchorbase]
\node at (0,0) {\includegraphics[height=5cm]{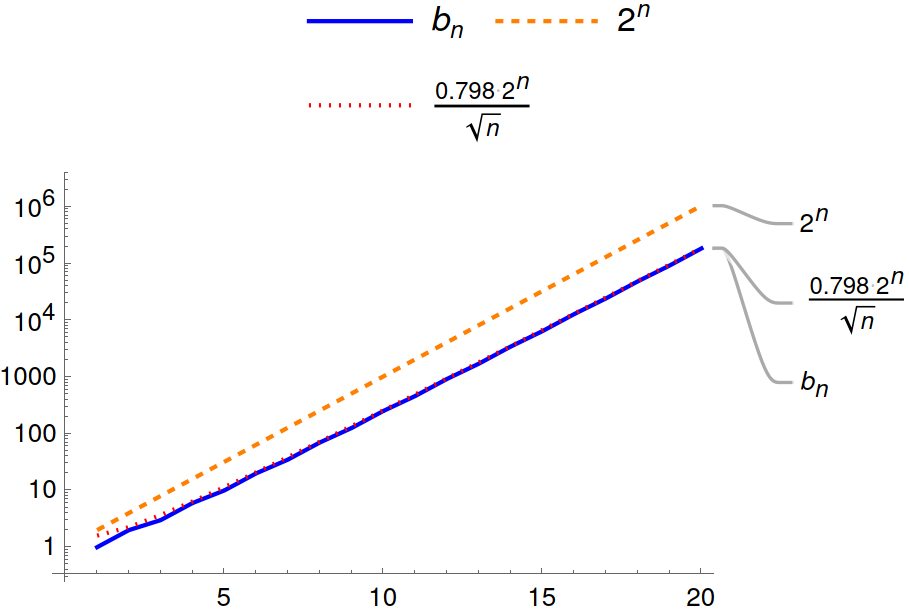}};
\end{tikzpicture}
,
\begin{tikzpicture}[anchorbase]
\node at (0,0) {\includegraphics[height=5cm]{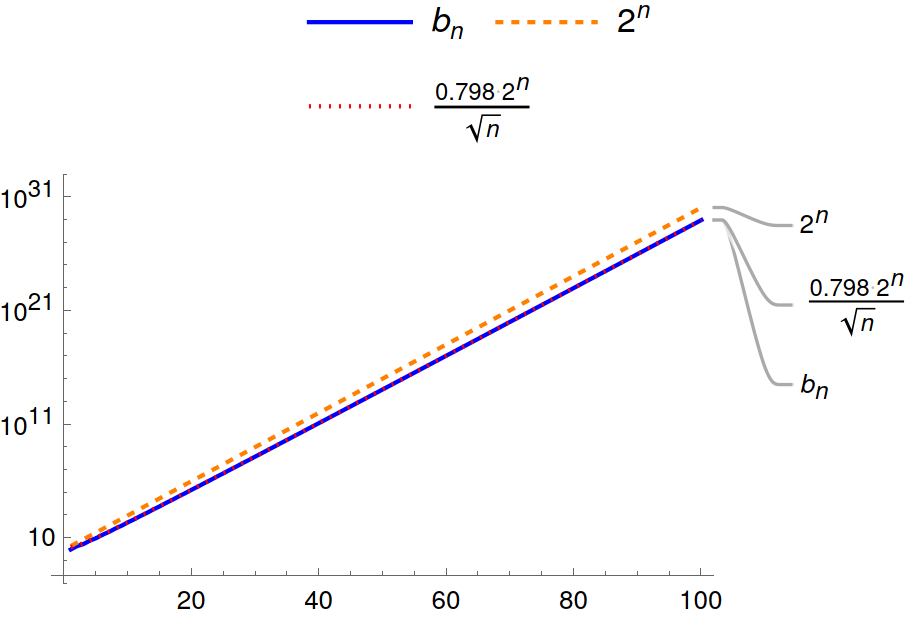}};
\end{tikzpicture}
\hspace*{-0.6cm}.
\end{gather*}
For a precise statement see \autoref{e:clt} below.
\end{Example}

A broader formulation of the ideas in \autoref{d:main} and \autoref{t:main} is the following setup. 

\begin{Notation}\label{n3:main}
Let $\cD$ be a $\bk$-linear Karoubian monoidal category that is Krull--Schmidt, with a $\bk$-linear faithful monoidal functor $F\colon\cD\to\cVect[\bK]$ to the category $\cVect[\bK]$ of finite 
dimensional vector spaces over a field extension $\bK$ of $\bk$.

{Note that, if $\bK$ is a finite extension of $\bk$, then the existence of $F$ implies that morphism spaces in $\cD$ are finite dimensional, so $\cD$ is automatically Krull--Schmidt.}
\end{Notation}

For any object $X\in\cD$, we can define $b_{n}^{X}$ similarly as in \autoref{d:main} as the number of indecomposable direct summands in $X^{\otimes n}$ and we have
\begin{gather*}
\beta^{\cD,X}:=\lim_{n\to\infty}\sqrt[n]{b_{n}^{X}}\leq\dim F(X).
\end{gather*}
Again, $\beta^{\cD,X}$ is well-defined by (a version of) Fekete's Subadditive Lemma.
We prove the following result, which also generalizes most of the examples in \autoref{t2:main}.

\begin{Theorem}\label{t3:main}
\leavevmode
\begin{enumerate}

\item If $\cD$ and $F$ are symmetric monoidal, then \autoref{t:main} holds, i.e.:
\begin{gather*}
\beta^{\cD,X}=\dim F(X).
\end{gather*}

\item Assume that $\mathrm{char}(\bk)=0$. If $\cD$ is symmetric and $F$ can be lifted to a symmetric monoidal functor $F^{\prime}\colon\cD\to\sVect[\bK]$ to the category of super vector spaces $\sVect[\bK]$, then \autoref{t:main} holds, i.e.:
\begin{gather*}
\beta^{\cD,X}=\dim F(X).
\end{gather*}

\item Let $\bk=\mC$ and $\cD=\Rep_{\mC}SL_{2}$. For every $m\in\mZ_{\geq 2}$, there exists a faithful monoidal functor $F_{m}\colon\cD\to\cVect[\mC]$, which sends the vector $SL_{2}$-representation $X$ to $\mC^{m}$. Hence, for $m\geq 3$ we have
\begin{gather*}
\beta^{\cD,X}=2<m=\dim F_{m}(X).
\end{gather*}

\end{enumerate}
\end{Theorem}

Here and throughout, $\Rep_{\bk}\Gamma$ denotes the category 
of finite dimensional (rational) $\Gamma$-representations over $\bk$.

\begin{Remark}\label{r3:main}
\leavevmode
\begin{enumerate}

\item Note that the functors $F_{m}$ in 
\autoref{t3:main}.(c) are not symmetric for $m\geq 3$, so 
\autoref{t3:main}.(a) does not apply.

\item Faithfulness of $F$ in \autoref{t3:main}.(a) is required to ensure the estimate
$b_{n}^{X}\leq\big(\dim F(X)\big)^{n}$ and cannot be dropped, since it is easy to construct counterexamples with Deligne's categories of \cite{De-cat-st}. Of course, even without faithfulness, the bound $\dim F(X)\le \beta^{\mathcal{D},X}$ remains valid.
\end{enumerate}
\end{Remark}

%%%%%%%%%%%%%%%%
% Acknowledgments
%%%%%%%%%%%%%%%%

\noindent\textbf{Acknowledgments.}
Want to thank Pavel Etingof for comments on a draft of this paper, in particular for \autoref{r:bounds} which is Pavel's observation, Andrew Mathas for help with the literature on symmetric groups and Schur algebras, Volodymyr Mazorchuk for email exchanges about monoids, and 
Jonathan Gruber and Arun Ram for discussions about growth rates.
{We also thank the referee for very valuable comments and remarks.}
We also thank the MFO workshop 2235 
``Character Theory and Categorification'' 
for bringing us together in Oberwolfach 
in August/September 2022 -- this project started during 
this fantastic workshop.

V.O. and D.T. were supported by their depressions.

%%%%%%%%%%%%%%%%%%%%%%%%%%%%%%%%%%%%%%%%%%%%%%

\section{The general linear group and consequences}\label{s:groups}

%%%%%%%%%%%%%%%%%%%%%%%%%%%%%%%%%%%%%%%%%%%%%%

Let us start by defining some of our main players:

\begin{Definition}\label{d:schemes}
We define an affine semigroup scheme (over $\bk$) to be a semigroup object in the category of affine $\bk$-schemes (the opposite of the category of commutative $\bk$-algebras). Equivalently, we can think of it as a representable functor from the category of $\bk$-algebras to the category of semigroups. Concretely, an affine semigroup scheme corresponds to a {commutative} bialgebra, potentially without counit. The special case of an affine monoid scheme corresponds precisely to a bialgebra with a counit, and the further special case of an affine group scheme corresponds to a Hopf algebra.

We have the same notions in the `super' version, by replacing the category of commutative $\bk$-algebras with the category of graded commutative $\BZ/2\BZ$-graded $\bk$-algebras.
\end{Definition}

Let $M\in\BZ_{>0}$, consider 
$GL_{M}$ and let $V_{M}$ be the 
tautological representation of $GL_{M}$ with $\dim V_{M}=M$.
Set $b_{n}^{M}=b_{n}^{GL_{M},V_{M}}$ and $\beta^{M}=\beta^{GL_{M},V_{M}}=\lim_{n\to\infty}\sqrt[n]{b_{n}^{M}}$.

\begin{Proposition}\label{p:glm}
For any $M\in\BZ_{>0}$ we have 
$\beta^{M}=M$.
\end{Proposition}

\begin{proof}
The case $M=1$ is immediate. In 
\autoref{ss:proofgltwo} 
and \autoref{ss:proofglm} 
we prove the statement for $M>1$.
\end{proof}

We also fix some notation:

\begin{Notation}\label{n:groundfield}
\leavevmode
\begin{enumerate}

\item Recall that $\bk$ denotes an arbitrary field, and we thus have that 
its characteristic $\mathrm{char}(\bk)\in\mN$. We always let $p\in \mZ_{>1}\cup\{\infty\}$ denote the additive order of $1\in\bk^{\times}$.
Thus, $\mathrm{char}(\bk)=p$ except that for $\mathrm{char}(\bk)=0$ we set $p=\infty$. When we state $\mathrm{char}(\bk)=p>0$ it thus unambiguously means positive characteristic.

\item Throughout this section, we will set $\Gamma_{M}=SL_{M}$.
Since the decomposition of $V^{\otimes n}_{M}$ into indecomposable summands is identical for $SL_M$ and $GL_M$, 
{we can focus just on $\Gamma_M$ in \autoref{p:glm}.}

\end{enumerate}
\end{Notation}

\begin{Remark}\label{r:tilting}
The proof of \autoref{p:glm} will crucially exploit the 
theory of tilting representations, see for example \cite[Part II.E]{Ja-reps-algebraic-groups} 
(note that the relevant section in \cite{Ja-reps-algebraic-groups} works over an arbitrary field 
of characteristic $\mathrm{char}(\bk)=p>0$) or \cite{AnStTu-cellular-tilting} and 
the appendix {of} its arXiv version for background. 
The point is that $V_{M}$ is tilting, 
and by abstract theory direct summands of tensor products of 
tilting representations are tilting. Thus, the summands 
of $V^{\otimes n}_{M}$ are tilting.
\end{Remark}

Before proving \autoref{p:glm}, we extract some of its consequences.

\begin{proof}[Proof of \autoref{t3:main}.(a)]
Let $\cD$ and $F$ be as in \autoref{t3:main}.(a), and let $X$ be an object of $\cD$. We have $\beta^{X}\leq\dim F(X)$ by the analog of \autoref{l:main}. We will show that
$b_{n}^{\dim F(X)}\leq b_{n}^{X}$,
which implies the claim by \autoref{p:glm}.

To this end, recall that $p=\mathrm{char}(\bk)\in\BN\cup\{\infty\}$.
Recall that simple representations of the symmetric group $S_{n}$ are labeled by $p$-regular partitions $\lambda$ of $n$, see e.g. \cite[Theorem 3.43]{Ma-hecke-schur} 
for an even more general statement. 
Let $D_{\lambda}$ be the simple
representation labeled by $\lambda$. {By Schur--Weyl duality, see e.g. \cite[E.17]{Ja-reps-algebraic-groups}, we have an epimorphism
\begin{gather*}
\bk S_n\twoheadrightarrow\End_{GL_M}(V_M^{\otimes n}).
\end{gather*}
Hence the number of indecomposable summands in $V_M^{\otimes n}$ equals the number of primitive idempotents of $\bk S_n$, in a fixed decomposition of the identity, that are not sent to zero. Such a decomposition always consists of $\dim D_\lambda$ idempotents corresponding to $D_\lambda$, for each $p$-regular partition $\lambda$.
Again by Schur-Weyl duality, this yields
\begin{gather*}
{b_{n}^{M}=\sum_{\lambda}\dim D_{\lambda},}
\end{gather*}
{where the sum runs over all $p$-regular partitions of $n$ with $\leq M$ rows.}}
By assumption, we have algebra morphisms
\begin{gather*}
\bk S_{n}\;\to\;\End_{\cD}(X^{\otimes n})\;\to\; 
\End_{\bk}\big(F(X)^{\otimes n}\big),
\end{gather*}
where the composite is the usual permutation action of $S_{n}$.
For any $p$-regular partition $\lambda$ let $e_{\lambda}\in\bk[S_{n}]$ be a primitive idempotent such that
$\bk[S_{n}]e_{\lambda}$ is the projective cover of $D_{\lambda}$ and let $\Pi_\lambda X$ be the direct summand $e_\lambda(X^{\otimes n})$ of $X^{\otimes n}$. Clearly $\Pi_\lambda X$ is not zero whenever $\Pi_{\lambda}F(X):=e_{\lambda}\big(F(X)^{\otimes n}\big)$ is not zero. 
Hence $X^{\otimes n}$ decomposes as a direct sum of $\Pi_{\lambda}X$ (where the latter need not be indecomposable)
with multiplicities $\dim D_{\lambda}$. By the above, 
this shows that
$b_{n}^{X}$ is indeed bounded below by $b_{n}^{\dim F(X)}$.
\end{proof}

\begin{proof}[Proof of \autoref{t2:main}.(a) -- affine group schemes]
For any affine group scheme (or any abstract group) $\Gamma$ we have
\begin{gather}\label{eq:estimate-tmain}
b_{n}^{\dim V}\leq
b_{n}^{\Gamma,V},
\end{gather}
since $V^{\otimes n}$, considered as a $\Gamma$-representation, is the restriction of the tensor power of the tautological $GL(V)$-representation under the usual homomorphism $\Gamma\to GL(V)$. 
Thus, in order to prove \autoref{t:main} we just need to combine the estimate $b_{n}^{\Gamma,V}\leq(\dim V)^{n}$ 
from \autoref{l:main} with \autoref{p:glm} and 
\autoref{eq:estimate-tmain}.
\end{proof}

\begin{proof}[Proof of \autoref{t2:main}.(a) -- affine semigroup schemes]
For any affine semigroup scheme (or any semigroup) $\Gamma$ let $\cD=\Rep_{\bk}\Gamma$. We have
a $\bk$-linear faithful monoidal functor $F\colon\cD\to\cVect[\bk]$ sending 
$X$ to its underlying $\bk$-vector space. Since this functor is 
symmetric, \autoref{t3:main}.(a) applies and we are done.
\end{proof}

There is also a proof for semigroups which does not rely on Schur--Weyl duality. Details will be given in \autoref{s:super} because Schur--Weyl duality fails in the super case in positive characteristic.

Finally, we conclude the list of consequences of \autoref{p:glm} with the case of groupoids. Consider a groupoid $(S:G)$ in the category of $\bk$-schemes, see \cite[\S 1.6]{Del90}, with source and target morphisms $G\rightrightarrows S$. If $G\to S\times S$ is faithfully flat, then we say $(S:G)$ is transitive. {A representation of a groupoid $(S:G)$ is a quasi-coherent sheaf on $S$ with an action of $G$.} Tensor products of representations are taken over $\cO_{S}$. We denote by 
$b_{n}^{(S:G),V}$ the analog of \autoref{d:main}.

\begin{Theorem}\label{CorGS}
Consider a transitive groupoid $(S:G)$ in $\bk$-schemes with a representation $V$ on a locally free $\cO_{S}$-module of finite rank. Then
\begin{gather*}
\pushQED{\qed} 
\beta^{(S:G),V}=\lim_{n\to\infty}\sqrt[n]{b_{n}^{(S:G),V}}=\mathrm{rank}(V).\qedhere
\popQED
\end{gather*}
\end{Theorem}

\begin{proof}
As observed in \cite[\S 1]{Del90}, the monoidal category $\Rep_{\bk}(S:G)$ of representations on locally free modules is actually abelian, and exact monoidal functors out of it are automatically faithful. For any field extension $\bK$ of $\bk$ for which $S(\bK)\not=0$, taking stalks at the corresponding point of $S$ therefore yields a faithful symmetric monoidal functor
\begin{gather*}
\Rep_{\bk}(S:G)\to\cVect[\bK],
\end{gather*}
which sends $V$ to a vector space of dimension $\mathrm{rank}(V)$. We can thus apply \autoref{t3:main}(a).
\end{proof}

%%%%%%%%%%%%%%%%%%%%%%%%%%%%%%%%%%%%%%%%%%%%%%

\subsection{Finite groups}\label{ss:finitegroups}

%%%%%%%%%%%%%%%%%%%%%%%%%%%%%%%%%%%%%%%%%%%%%%

\begin{Notation}\label{n:landaubachmann}
Recall the {Bachmann--Landau} notation, which we adjust as follows. 

A function $f$ satisfies $f\in\Theta^{\prime}(g)$ if there exists a constant $A\in\BR_{>0}$ such that $A\cdot g(n)\leq f(n)\leq g(n)$ for 
all $n_{0}<n$ for some fixed $n_{0}\in\BN$. Without the prime this is the classical {Bachmann--Landau} notation where $A\cdot g(n)\leq f(n)\leq B\cdot g(n)$ for $A,B\in\BR_{>0}$.
Similarly, but only for either the lower or the upper bound, 
we write $f\in\Omega(g)$ and $f\in O^{\prime}(g)$, respectively. Finally, we write $f\in \Omega^{\prime}(g)$ if  $g\in O'(f)$.

We also use $f\sim g$, $f$ is asymptotic to $g$, meaning 
$\lim_{n\to\infty}f(n)/g(n)=1$ (with $g(n)\neq 0$ for all $n\gg 1$). 

\end{Notation}

\begin{Proposition}\label{p:finitegroup}
If $\Gamma$ is a finite group, then
\begin{gather*}
b_{n}^{\Gamma,V}\in\Theta^{\prime}\big((\dim V)^{n}\big).
\end{gather*}
\end{Proposition}

\begin{proof}
Without loss of generality, we assume that $\Gamma$ acts faithfully on $V$. By
\cite[Theorem 1]{BrKo-tensor-group}, there exists $r\in\BN$ for which $V^{\otimes r}$ contains a projective direct summand. Using that projective representations 
form a tensor ideal, and the fact that projective indecomposables $P$ satisfy $\dim P\leq|\Gamma|$, it then follows that
\begin{gather*}
\frac{(\dim V)^{n}}{(\dim V)^{r}|\Gamma|}
\leq b_{n}^{\Gamma,V},
\end{gather*}
which, by \autoref{l:main}, concludes the proof.
\end{proof}

This implies \autoref{t2:main} for finite groups.

\begin{Example}\label{e:projective}
Let $\Gamma_{2}=SL_{2}(\bk)$ and let $V_{2}$ be its vector representation.
For $\bk=\BC$ we have seen in \autoref{e:main} that 
$b_{n}^{\Gamma_{2},V_{2}}\sim A\cdot 2^{n}/\sqrt{n}$ for $A\in\BR_{>0}$.
Let $\BF_{p^{r}}$ be the finite Galois field with $p^{r}$ elements. 
For $\bk=\overline{\BF_{p^{r}}}$, since the number of indecomposable summands is bounded from above by the number of indecomposable summands 
for $\bk=\BC$, 
\autoref{e:main} implies that the growth 
rate of $b_{n}^{\Gamma_{2},V_{2}}$ is bounded by $A\cdot 2^{n}/\sqrt{n}$
for $A\in\BR_{>0}$. With contrast, in the case 
$\bk=\BF_{p^{r}}$ \autoref{p:finitegroup} shows that 
$b_{n}^{\Gamma_{2},V_{2}}\in\Theta^{\prime}(2^{n})$ and hence, we do not have an upper bound by $A\cdot 2^{n}/\sqrt{n}$
for $A\in\BR_{>0}$.
This in turn implies that 
Schur--Weyl duality fails over finite fields. 
The latter was observed in \cite[Section 5]{BeDo-finite-schur-weyl}.
\end{Example}

\begin{Remark}\label{r:projective}
The argument in the proof of \autoref{p:finitegroup} is that, for $n\gg 0$, 
most indecomposable summands are projective, their number of appearance
can be bound from below and grows already fast enough.
As we will see in \autoref{ss:proofglm} below, for the special and general linear groups the role of projective representations will be played 
by certain tilting representations.
\end{Remark}

\begin{Remark}\label{r:monoid}
Let $\Gamma$ be a monoid or semigroup. Note that $\bk\Gamma$ is not a Hopf algebra, but only a bialgebra 
(potentially without unit). In particular, we can tensor representations, 
but the projective representation do not form 
a tensor ideal, see \cite[Exercise 17.15]{St-rep-monoid} for an explicit counterexample.
Thus, the arguments in the proof of 
\autoref{p:finitegroup} do not apply, not even 
for the case of finite monoids or semigroups. From this point of view it is surprising that \autoref{t:main} remains valid.
\end{Remark}

%%%%%%%%%%%%%%%%%%%%%%%%%%%%%%%%%%%%%%%%%%%%%%

\subsection{Proof of \autoref{p:glm} -- semisimple case}\label{ss:semisimple}

%%%%%%%%%%%%%%%%%%%%%%%%%%%%%%%%%%%%%%%%%%%%%%

Let $\mathrm{char}(\bk)=0$, which we call the semisimple case.

\begin{Notation}\label{r:weyl}
For $m_{1},\dots,m_{M-1}\in\BN^{M-1}$ we denote by 
$\Delta(m_{1},\dots,m_{M-1})$ the Weyl representations of 
$\Gamma_{M}=SL_{M}$ of highest weight 
$(m_{1},\dots,m_{M-1})$. This is in terms of the fundamental weights 
meaning that the highest weight is $\sum_{i=1}^{M-1}m_{i}\omega_{i}$ where the $\omega_{i}$ the fundamental weights.

These are $\Gamma_{M}$-representations
defined integrally and these are simple for $\mathrm{char}(\bk)=0$, and we 
also have $\Delta(1,0,\dots,0)=V_{M}$.
See \cite[Part II]{Ja-reps-algebraic-groups} for some background.
\end{Notation}

For $\mathrm{char}(\bk)=0$ the tensor product $V^{\otimes n}_{M}$ 
decomposes into the simple summands for $m_{1},\dots,m_{M-1}\in\BN^{M}$. 
Recall Weyl's character formula, see e.g.  
\cite[Section 24]{FuHa-representation-theory}, which shows that 
\begin{gather*}
\dim\Delta(m_{1},\dots,m_{M-1})\in\BN`[m_{1},\dots,m_{M-1}],
\end{gather*}
i.e. $\dim\Delta(m_{1},\dots,m_{M-1})$ is a polynomial in $m_{1},\dots,m_{M-1}$.

\begin{Example}\label{e:weyl}
For example, for $M=3$ one has
$\dim\Delta(m_{1},m_{2})=\frac{1}{2}(m_{1}+1)(m_{2}+1)(m_{1}+m_{2}+2)$.
\end{Example}

The following implies \autoref{p:glm} for $\mathrm{char}(\bk)=0$:

\begin{Proposition}\label{p:semisimple}
We have
\begin{gather*}
{\Omega\big(M^{n}/n^{M(M-1)/2}\big)}\ni b_{n}^{M}\in O^{\prime}(M^{n}).
\end{gather*}
\end{Proposition}

\begin{proof}
All the weights of the representation $V^{\otimes n}_{M}$
are bounded by $n$ in the sense that any coefficient of 
the expansion with respect to fundamental weights
is less than $n$ in absolute value, meaning that $m_{1}+\dots+m_{M-1}\leq n$.
Thus, all simple
summands of $V^{\otimes n}_{M}$ have dimensions bounded by a polynomial in $n$, and a closer look at Weyl's character formula then implies that the polynomial is of degree $\frac{M(M-1)}{2}$. See \autoref{e:weyl} 
for an example. Thus, from this and \autoref{l:main} we get 
\begin{gather*}
\frac{M^{n}}{f(n)}\leq b_{n}^{M}\leq M^{n},
\quad f(x)\in\BN[x],\;\deg f=\tfrac{M(M-1)}{2}.
\end{gather*}
The claim follows.
\end{proof}

\begin{Example}\label{e:clt}
We now strengthen \autoref{p:semisimple} for $\Gamma_{2}$. For a $\Gamma_{2}$-representation $W$, with weight spaces $\{W_{i}\subset W | i\in\mZ\}$, its character is the Laurent polynomial with non-negative coefficients $\mathrm{ch}\,W=\sum_{i}(\dim W_{i})v^{i}\in\mN[v,v^{-1}]$.	

Let $V$ be a representation of $\Gamma_{2}$ over $\bk$ with $\mathrm{char}(\bk)=0$. Then, by the classification of simple $SL_{2}(\bk)$-representations, 
$b_{n}^{G,V}$ equals the sum of dimensions of zero weight space and one weight space of $V^{\otimes n}$. For example if $V=V_{2}$ is the vector representation,
then $\mathrm{ch}\,V_{2}=v+v^{-1}$ and $\mathrm{ch}\,V_{2}^{\otimes n}=(v+v^{-1})^{n}$.
This implies that $b_{n}^{\Gamma_{2},V_{2}}$ is the 
constant term or the coefficient of $v$ of $(v+v^{-1})^{n}$, 
depending on parity of $n$.
Using the binomial theorem we see that
\begin{gather*}
b_{n}^{\Gamma_{2},V_{2}}=
\begin{cases}
\binom{n}{n/2}& n\;\text{even},
\\
\binom{n}{(n-1)/2}& n\;\text{odd}.
\\
\end{cases}
\end{gather*}
By applying Stirling's formula we get that asymptotically
\begin{gather*}
b_{n}^{\Gamma_{2},V_{2}}\sim\sqrt{2/\pi}\cdot\frac{2^{n}}{\sqrt{n}},
\end{gather*}
which implies the claim in \autoref{e:main}.
Note that this is better than 
what we get from \autoref{p:semisimple} for $M=2$. 
That is, the lower bound $2^{n}/n$ is what we get from \autoref{p:semisimple} and
Mathematica's log plot gives:
\begin{gather*}
\begin{tikzpicture}[anchorbase]
\node at (0,0) {\includegraphics[height=5cm]{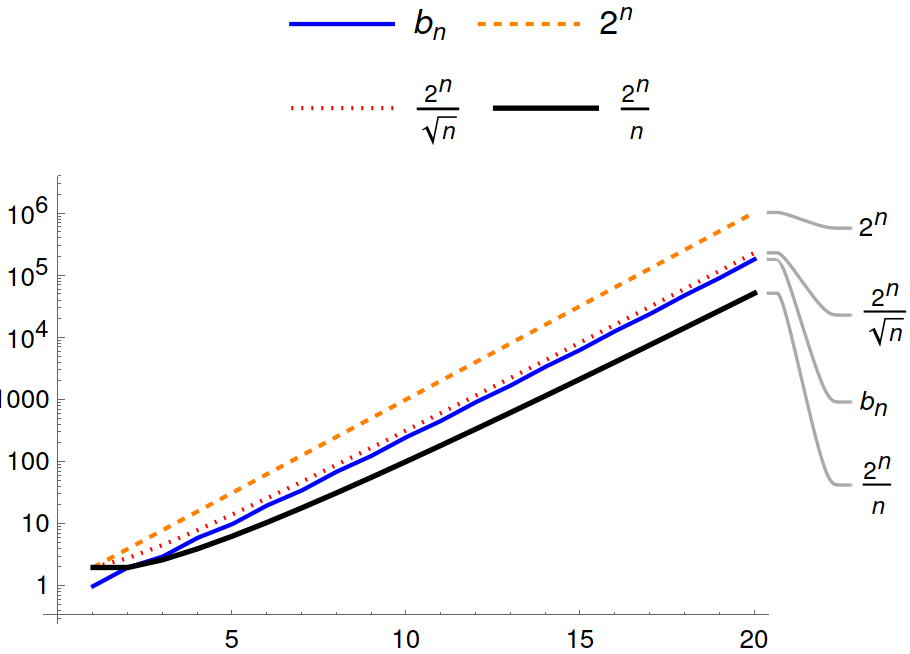}};
\end{tikzpicture}
,
\begin{tikzpicture}[anchorbase]
\node at (0,0) {\includegraphics[height=5cm]{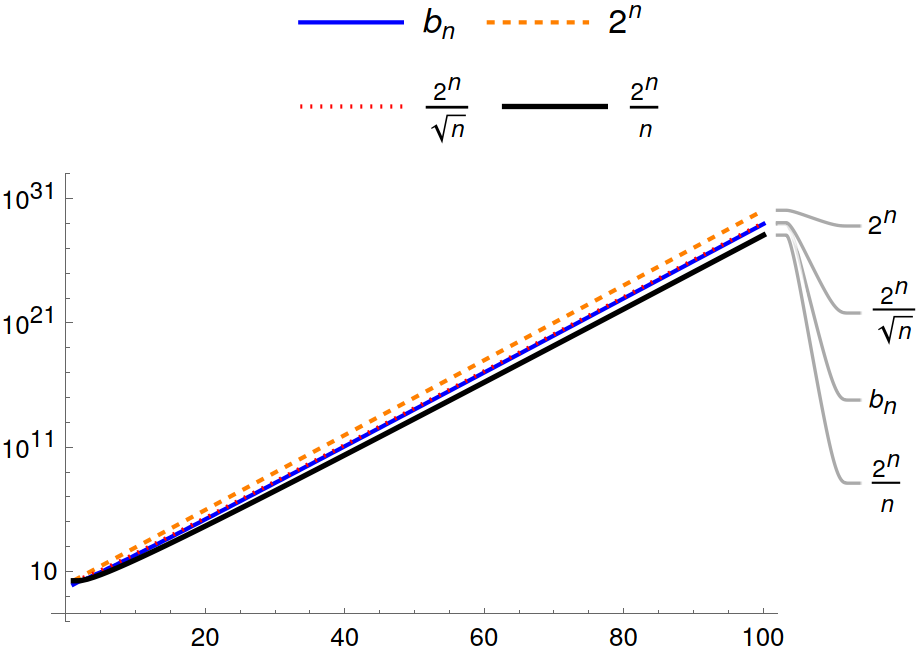}};
\end{tikzpicture}
.
\end{gather*}

More generally, let $V$ be any nontrivial representation of 
$\Gamma_{2}$. Then 
$b_{n}^{\Gamma_{2},V}$ equals the sum of the constant term and the coefficient of $v$ in $(\mathrm{ch}\,V)^{n}$. The asymptotic 
of this number can be computed by using the central limit theorem.
We get that
\begin{gather*}
b_{n}^{\Gamma_{2},V}\sim
(\dim V)^{n}/(2\pi nA^{2})^{1/2},
\end{gather*}
where $A=A(V)\in\BR_{>0}$ is an easily computable constant depending on $V$.
In the case when $V$ is simple multinomials appear and one can use e.g. the results from \cite{Eg-stirling-binomial}. The same approach applies in general.
\end{Example}

\begin{Remark}\label{r:manyrelatedresults}
In the semisimple case many related results are known, in particular for Lie algebras and Lie groups, see e.g. 
\cite{PoRe-mult-large-tensor-powers} for a recent publication. That paper 
studies the problem of finding the asymptotic of multiplicities of fixed simple representations instead of all simple representations. 
\end{Remark}

A growth rate $b_{n}^{\Gamma,V}\in\Theta^{\prime}\big((\dim V)^{n}\big)$ as in \autoref{p:finitegroup} is very rarely the case, as the following result indicates:

\begin{Proposition}\label{p:fastgrowth}
Recall that $\mathrm{char}(\bk)=0$. For an abstract group $\Gamma$ with a finite dimensional representation $V$, the following are equivalent:
\begin{enumerate}

\item We have $b_{n}^{\Gamma,V}\in\Theta^{\prime}\big((\dim V)^{n}\big)$.

\item The connected component of the Zariski closure of the image of $\Gamma$ in $GL(V)$ is a torus.

\end{enumerate}
\end{Proposition}

\begin{proof}
We start by proving that (a) implies (b). Note that, with $W=\bar{\bk}\otimes_{\bk}V$ equipped with the canonical structure of a $\Gamma$-representation over $\bar{\bk}$, we have $b_{n}^{\Gamma,V}\leq b_{n}^{\Gamma,W}$. By definition, an algebraic group is a torus over $\bk$ if and only if its extension of scalars to $\bar{\bk}$ is a torus (a finite product of copies of the multiplicative group). Consequently, for this implication, we might as well assume that $\bk$ is algebraically closed.

Replacing $\Gamma$ by the Zariski closure of its image in $GL(V)$ does not change the
numbers $b_{n}^{\Gamma,V}$, so we can assume that $\Gamma$ is 
an algebraic group and $V$ is a faithful
$\Gamma$-representation. 

Now we argue by contradiction. We assume that (a) is satisfied. If (b) is not satisfied, then by \cite[Corollary 17.25]{Milne} and the fact that every connected one-dimensional unipotent algebraic group is isomorphic to the additive group, it follows that $\Gamma$ contains a copy of the additive group $\BG_{a}$. 
We can restrict $V$ to $\BG_{a}$ and get
\begin{gather*}
A\cdot(\dim V)^{n}\leq b_{n}^{\Gamma,V}\leq b_{n}^{\BG_{a},V},
\end{gather*}
for some $A\in\BR_{>0}$. We can apply the Jacobson--Morozov theorem and find $\Gamma_{2}=SL_{2}\subset GL(V)$ containing $\BG_{a}$ and such that
$b_{n}^{\BG_{a},V}=b_{n}^{\Gamma_{2},V}$, and thus, 
\begin{gather*}
A\cdot(\dim V)^{n}\leq b_{n}^{\Gamma_{2},V},
\end{gather*}
for the same $A\in\BR_{>0}$.
This gives a contradiction with \autoref{e:clt}, concluding the proof of this direction.

Now we prove that (b) implies (a). For a finite field extension $\bK$ of $\bk$, we can again consider a $\Gamma$-representation on $U=\bK\otimes_{\bk}V$. 
Restricting the $\bK$-vector space $U^{\otimes_{\bK} n}$  to $\bk$ yields a direct sum of $[\bK:\bk]$ copies of $V^{\otimes n}$. Since this isomorphism respects the $\Gamma$-action, it follows that
\begin{gather*}
b_{n}^{\Gamma,U}\leq[\bK:\bk]\cdot b_{n}^{\Gamma, V}.
\end{gather*}
We can again replace $\Gamma$ by the Zariski closure of its image. Since any torus splits after a finite field extension, we can thus assume that
\begin{gather*}
\Gamma=(\BG_{m}^{\times d})\rtimes H	
\end{gather*}
for $d\in\mN$ and a finite group $H$. It is well-known that for such groups the (rational) representation theory is semisimple and the dimension of the simple representations are bounded by $|H|$. Conclusion (a) from this 
in the same way as in \autoref{p:finitegroup}.
\end{proof}

%%%%%%%%%%%%%%%%%%%%%%%%%%%%%%%%%%%%%%%%%%%%%%

\subsection{Proof of \autoref{p:glm} -- for \texorpdfstring{$M=2$}{M=2}}\label{ss:proofgltwo}

%%%%%%%%%%%%%%%%%%%%%%%%%%%%%%%%%%%%%%%%%%%%%%

The case $M=2$ is special since we have full access 
to the characters of tilting representations
and these are the direct summands of $V^{\otimes n}_{2}$. 
As before, $\Gamma_{2}=SL_{2}$.

It is crucial that $V^{\otimes n}_{2}$ is 
a tilting $\Gamma_{2}$-representation, see \autoref{r:tilting}. Thus, its direct summands
are indecomposable tilting representations $T(m)$ parameterized by dominant weights $m\in\BN$.
Moreover, the $\Gamma_{2}$-representation $V^{\otimes n}_{2}$ decomposes into direct summands $T(m)$ with $m\leq n$. 

\begin{Example}\label{e:sl2}
Donkin's tensor product theorem \cite[Proposition 2.1]{Do-tilting-alg-groups} allows us to describe tilting characters explicitly for $\Gamma_2$. As observed in e.g. \cite{TuWe-quiver-tilting} or 
\cite[Section 2]{SuTuWeZh-mixed-tilting}, we can reformulate 
Donkin's result for $\Gamma_{2}$ as follows.
Let $m+1=a_{d}p^{d}+\cdots+a_{1}p+a_{0}=(a_{d},\dots,a_{0})$ be the $p$-adic expansion of $m+1$ where $a_{i}\in\{0,\dots,p-1\}$ and $a_{d}\neq 0$. (The $a_{0}$ digit is the one for $p^{0}$, where we use the convention that $\infty^{0}=1$.)
We use the convention for characters from \autoref{e:clt}. For $b\in\BN$ write $[b]_{x}=x^{-(b-1)}+x^{-(b-3)}+\dots+x^{b-3}+x^{b-1}$.
We then have
\begin{gather*}
\mathrm{ch}\,T(m)=
[a_{d}p^{d}]_{v}\cdot\prod_{a_{i}\neq 0,i\neq d}
[2]_{v^{a_{i}p^{i}}}.
\end{gather*}
For example, for $m=52$ and $p=2$ we have 
$m+1=(1,1,0,1,0,1)$ so
\begin{gather*}
\mathrm{ch}\,T(52)=
[p^{5}]_{v}\cdot 
[2]_{v^{2^{4}}}[2]_{v^{2^{2}}}[2]_{v^{2^{0}}}.
\end{gather*}
In particular, for $v=1$ we get $\dim T(52)=256$.
\end{Example}

This example implies:

\begin{Lemma}\label{l:sl2} 
Let $\alpha=1+(\log_{2}p)^{-1}$. Then we have
\begin{gather*}
\dim T(m)\leq(m+1)^{\alpha}.
\end{gather*}
\end{Lemma}

\begin{proof} 
It follows from \autoref{e:sl2} that
\begin{gather*}
\dim T(m)=2^{k}a_{d}p^{d},
\end{gather*}
where $k$ is the number of non-zero digits among the $a_{d-1},\dots,a_{1},a_{0}$ in the $p$-adic extension of $m+1$. This implies
$\dim T(m)\leq(m+1)^{\alpha}$: First, we have $a_{d}p^{d}\leq m+1$ so it remains to argue that $2^{k}\leq(m+1)^{\log_{2} 2/\log_{2} p}$.
Note secondly that $2^{k}\leq 2^{d-1}$ and $2^{d-1}=(m+1)^{b}$ 
for $b=(d-1)\log_{2} 2/\log_{2}(m+1)$. However, $m+1=a_{d}p^{d}+\cdots+a_{1}p+a_{0}$ 
for $a_{d}\neq 0$ which gives $\log_{2}(m+1)\geq\log_{2}(p^{d-1})=(d-1)\log_{2} p$ so that $b\leq\log_{2} 2/\log_{2} p$. The result follows.
\end{proof}

\begin{Remark}\label{r:sl2}
The number $\alpha=1+(\log_{2} p)^{-1}$ in \autoref{l:sl2} 
converges to $1$ for $p\to\infty$, and $p=\infty$ is the semisimple
case where $\dim T(m)=m+1$.
\end{Remark}

\begin{Example}\label{r:sl2tdim}
For $p=2$ (left) and $p=3$ (right) we get
\begin{gather*}
\begin{tikzpicture}[anchorbase]
\node at (0,0) {\includegraphics[height=4.86cm]{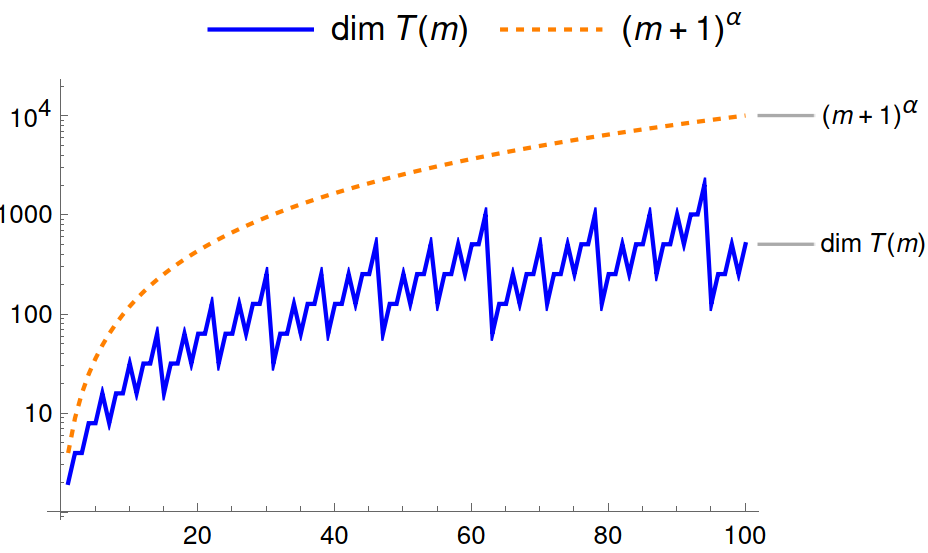}};
\end{tikzpicture}
,\hspace*{-0.19cm}
\begin{tikzpicture}[anchorbase]
\node at (0,0) {\includegraphics[height=4.86cm]{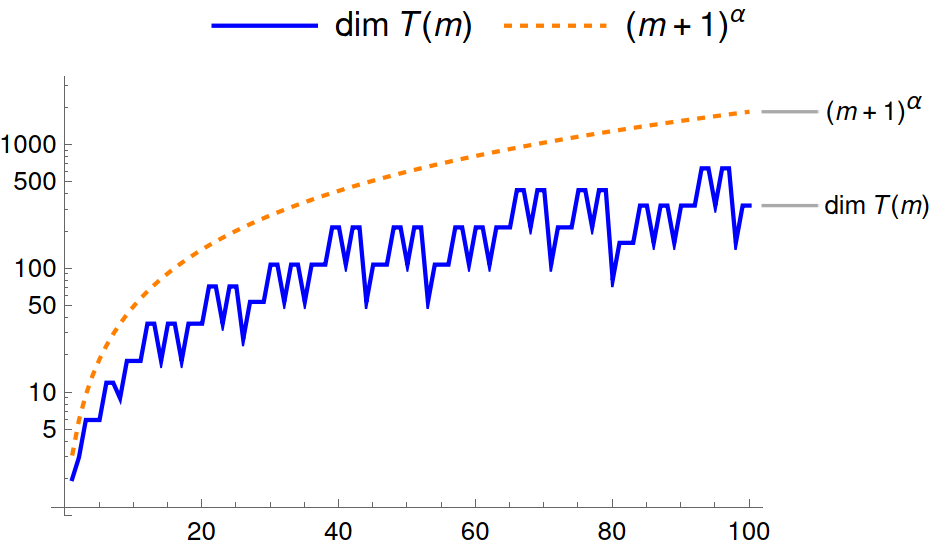}};
\end{tikzpicture}
\end{gather*}
which are again Mathematica log plots.
\end{Example}

The following is a finer result than \autoref{p:glm} itself, and thus, also implies \autoref{p:glm}.

\begin{Proposition}\label{p:sltwo}
We have
\begin{gather*}
\Omega^{\prime}\big(2^{n}/n^{2}\big)\ni b_{n}^{2}\in O^{\prime}(2^{n}).
\end{gather*}
\end{Proposition}

\begin{proof}
Since $1\leq\alpha\leq 2$, \autoref{l:main} and \autoref{l:sl2} give us
\begin{gather*}
2^{n}/(n+1)^{2}\leq 2^{n}/(n+1)^{\alpha}\leq b_{n}^{2}\leq 2^{n},
\end{gather*}
and the statement follows.
\end{proof}

%%%%%%%%%%%%%%%%%%%%%%%%%%%%%%%%%%%%%%%%%%%%%%

\subsection{Proof of \autoref{p:glm} -- for \texorpdfstring{$M\geq 2$}{M>=2}}\label{ss:proofglm}

%%%%%%%%%%%%%%%%%%%%%%%%%%%%%%%%%%%%%%%%%%%%%%

We assume $M\geq 2$ and $\mathrm{char}(\bk)=p>0$ as the case 
$\mathrm{char}(\bk)=0$ is dealt with in \autoref{ss:semisimple}.
In particular, the results from \cite[Part II.E]{Ja-reps-algebraic-groups} apply.
Let $\Gamma_{M}=SL_{M}$.

\begin{Remark}\label{r:donkincone}
The bound as in \autoref{l:sl2} is unavailable for $M\geq 3$; the 
billiards conjecture in \cite{LuWi-billiard-conjecture} and \cite{Je-correction-tilting-sl3} suggests that
dimensions of tilting representations along the 
boundary grow exponentially already for 
$\Gamma_{3}=SL_{3}(\bk)$.

Our proof below ``ignores'' these tilting representations: 
we argue that we already have enough summands in the part where Donkin's tensor product formula applies up to a certain degree.
\end{Remark}

We will use the following to only consider the case when $M$ is odd 
since this case has slightly nicer combinatorics:

\begin{Lemma}\label{l:avoideven}
If \autoref{p:glm} holds for $M+1$, then it holds for $M$ as well.
\end{Lemma}

\begin{proof}
Recall that the $\Gamma_{M+1}$-representation $V_{M+1}$ restricts to
the $\Gamma_{M}$-representation $V_{M}\oplus\bk$
under the usual embedding $\Gamma_{M}\hookrightarrow\Gamma_{M+1}$. It follows that
\begin{gather*}
b_{n}^{M+1}\leq\sum_{i=0}^{n}\binom{n}{i}b_{i}^{M}.
\end{gather*}
Using that \autoref{p:glm} holds for $M+1$ we get:
\begin{gather*}
M+1\leq\lim_{n\to\infty}a_{n},\quad a_{n}=\sqrt[n]{\sum_{i=0}^{n}\binom{n}{i}b_{i}^{M}}.
\end{gather*}
We claim that thus $\lim_{n\to\infty}\sqrt[n]{b_{n}^{M}}=M$, as required.
This can be seen as follows. 
Assume that for some fixed $\epsilon\in\BR_{>0}$ and all $\delta_{N}\in\BR_{>0}$ there exists $N\in\BN$ such that
$|(M-\epsilon)-(b_{i}^{M})^{1/i}|<\delta_{N}$ for $i\geq N$.

Since $\sum_{i=0}^{n}\binom{n}{i}(M-\epsilon)^{i}=(M-\epsilon+1)^{n}$ 
and $b_{i}^{M}\approx(M-\epsilon)^{i}$ for $i\gg 0$, it follows that $a_{n}=M-\epsilon+1-\delta_{n}^{\prime}$ for some $\delta_{n}^{\prime}\in\BR_{>0}$ 
with $\lim_{n\to\infty}\delta_{n}^{\prime}=0$.
Hence, $\lim_{n\to\infty}a_{n}=M-\epsilon+1<M+1$, which contradicts
$M+1\leq\lim_{n\to\infty}a_{n}$ and the proof completes.
\end{proof}

Recall that the category of finite dimensional $\Gamma_{M}=SL_{M}$-representations,
considered as an abelian category, has a direct summand $ST_{r}(\Gamma_{M})=ST_{r}^{p}(\Gamma_{M})$ consisting of representations which are linked
with the Steinberg representation $St_{r}=St_{r}^{p}=T\big((p^{r}-1)\rho\big)$ 
(note that these depend on $p$), see \cite[Section 3.5]{An-steinberg-linkage}. These Steinberg representations 
are tilting and Weyl representations at the same time, and we will use this below. 

We plan to choose $r=r(n)$ in such a way that the
number of summands of $V^{\otimes n}_{M}$ from $ST_{r}(\Gamma_{M})$ is still about $M^{n}$.
Let us first estimate the number of occurrences of $St_{r}$ as a subquotient of a good filtration of $V^{\otimes n}_{M}$.
This number depends only on the character of $V^{\otimes n}_{M}$ and hence, is independent of the
characteristic in the sense that the characters of both, $V^{\otimes n}_{M}$ and $St_{r}\cong\Delta\big((p^{r}-1)\rho\big)$, are as in characteristic zero. 

In fundamental weight coordinates and $SL_{M}$ notation, we let $\rho_{SL_{M}}=(1,\dots,1)$ 
and there are choices involved how to lift this to $GL_{M}$ notation in standard coordinates. 
We will use
\begin{gather*}
\rho=(\rho_{1},\dots,\rho_{M})=
\bigg(\frac{M-1}{2},\frac{M-3}{2},\dots,\frac{-M+3}{2},\frac{-M+1}{2}\bigg).
\end{gather*}
Now the number of times that $St_{r}$ appears in $V_{M}^{\otimes r}$ over $SL_{M}$ is at least the number of times it appears when we work over $GL_{M}$, so we estimate the latter number.

In characteristic zero we can compute the involved characters via Schur--Weyl duality by applying the hook length formula to the partition
\begin{gather*}
\lambda=
\lambda(p,r)=
\frac{n}{M}(1,\dots,1)+(p^{r}-1)\rho,
\end{gather*}
where we from now on
assume that $n$ is divisible by $M$ (which is sufficient to calculate the limit of $\sqrt[n]{b_n^M}$), that $M$ is odd (which is justified by \autoref{l:avoideven}) and that $\frac{(p^{r}-1)(M-1)}{2}\leq\frac{n}{M}$. 

\begin{Example}\label{e:partition}
Let $n=12$, $M=3$, $p=2$ and $r=2$. Note that $3=\frac{(p^{r}-1)(M-1)}{2}\leq\frac{n}{M}=4$ is satisfied.
Then $\lambda$ is	
\begin{gather*}
\lambda=
\begin{ytableau}
~ & ~ & ~ & ~ & ~ & ~ & ~
\\
~ & ~ & ~ & ~
\\
~
\end{ytableau}
\,.
\end{gather*}
This is the partition $(7,4,1)$ so that the row differences are $p^{r}-1=3$.
\end{Example}

Now let us make the following concrete choice for $r$: 
\begin{gather*}
r(n)=\lfloor\log_{p}(\sqrt{n})\rfloor.
\end{gather*}

\begin{Remark}\label{r:choiceofr} 
In fact, we could use 
$r(n)=\big\lfloor\log_{p}\big(f(n)\big)\big\rfloor$ for every function $f$ which grows slower than $n$. The choice $r(n)=\lfloor\log_{p}(\sqrt{n})\rfloor$ is mostly for convenience 
as the formulas come out nicely.
\end{Remark}

The hook formula implies that, up to factors which will not contribute to the limit of the $n$th root, the number of times that $St_{r(n)}=\Delta(\lambda)$
for $\lambda$ as above appears in $V^{\otimes n}$ over $GL_{M}$ is approximately
\begin{gather*}
a(n)=\frac{n!}{\big(\frac{n}{M}+(p^{r(n)}-1)\rho_{1}\big)!\dots\big(\frac{n}{M}+(p^{r(n)}-1)\rho_{M}\big)!}.
\end{gather*}
Let us write $x_{i}=x_i(n)=\frac{n}{M}+(p^{r(n)}-1)\rho_{i}$ (so $\sum_{i}x_{i}=n$). To approximate the above formula 
recall that, for all $a\in\BZ_{\geq 1}$, we have
\begin{gather*}
\sqrt{2\pi a}\left(\frac{a}{e}\right)^{a}e^{\frac{1}{12a+1}} 
<a!<
\sqrt{2\pi a}\left(\frac{a}{e}\right)^{a}e^{\frac{1}{12a}}.
\end{gather*}
Hence, we get that
\begin{gather*}
\colorbox{blue!25}{\mystrut$\frac{\sqrt{2\pi n}\left(n\right)^{n}}
{\prod_{i}\sqrt{2\pi x_{i}}\left(x_{i}\right)^{x_{i}}}$}
\cdot
\colorbox{spinach!25}{\mystrut$e^{\frac{1}{12n+1}-\sum_{i}\frac{1}{12x_{i}}}$}
<a(n)<
\colorbox{blue!25}{\mystrut$\frac{\sqrt{2\pi n}\left(n\right)^{n}}
{\prod_{i}\sqrt{2\pi x_{i}}\left(x_{i}\right)^{x_{i}}}$}
\cdot
\colorbox{tomato!25}{\mystrut$e^{\frac{1}{12n}-\sum_{i}\frac{1}{12x_{i}+1}}$}
.
\end{gather*}
{We claim that 
\begin{gather*}
\colorbox{spinach!25}{\mystrut$e^{\frac{1}{12n+1}-\sum_{i}\frac{1}{12x_{i}}}$}
\quad\text{and}\quad
\colorbox{tomato!25}{\mystrut$e^{\frac{1}{12n}-\sum_{i}\frac{1}{12x_{i}+1}}$}
\end{gather*}
then converge to one if $n\to\infty$.} Indeed, because of our choice for $r(n)$, $|(p^{r(n)}-1)\rho_{i}|$ is bounded from above by $B\cdot\sqrt{n}$ for some $B\in\BR_{>0}$ and thus, 
$\lim_{n\to\infty}x_{i}=\infty$. We then can also see that 
the {exponents
\begin{gather*}
\frac{1}{12n+1}-\sum_{i}\frac{1}{12x_{i}}
\quad\text{and}\quad
\frac{1}{12n}-\sum_{i}\frac{1}{12x_{i}+1}
\end{gather*}
converge} to zero, and the claim follows.
Consequently, using also $x_{i}(n)\sim\frac{n}{M}$, we find
\begin{gather*}
a(n)\sim\frac{1}{(2\pi)^{\frac{M-1}{2}}}\cdot\frac{n^{n+1/2}}{\prod_{i} x_{i}^{x_{i}+1/2}}\sim \frac{M^{M/2}}{(2\pi)^{\frac{M-1}{2}}}\cdot \frac{n^{n+\frac{1-M}{2}}}{\prod_{i}x_{i}^{x_{i}}}.
\end{gather*}
Let $f(n)$ denote a function with 
$f(n)\in\Theta(n^{-1/2})$ in {Bachmann--Landau} notation, see \autoref{n:landaubachmann}.
We get
\begin{gather*}
a(n)
\sim\colorbox{blue!25}{\mystrut$A\cdot e^{f(n)}n^{(1-M)/2}$}\colorbox{spinach!25}{\mystrut$M^{n}$},
\end{gather*}
for $A\in\BR_{>0}$.
This can be seen by using that
\begin{gather*}
\prod_{i}1/x_{i}^{x_{i}}\sim B\cdot e^{g(n)}n^{-n}M^{n}
\end{gather*}
for some $B\in\BR_{>0}$ and $g(n)\in\Theta(n^{-1/2})$.

Thus, since the limit $n\to\infty$ of the $n$th root of the 
(marked in a blueish color) left-hand side is one, we see that $n$th root of this sequence converges to $M$ and we get:
\begin{gather*}
\lim_{n\to\infty}\sqrt[n]{a(n)}
=\lim_{n\to\infty}\sqrt[n]{A\cdot e^{f(n)}n^{(1-M)/2}M^{n}}
=M.
\end{gather*}
Now let $t_{n}$ be  the 
total dimension of summands 
of $V^{\otimes n}_{M}$ which are in $ST_{r(n)}$. Clearly, we have $a(n)\leq t_{n}\leq(\dim V)^{n}$ and thus in conclusion
\begin{gather*}
\lim_{n\to\infty}\sqrt[n]{t_{n}}=M.
\end{gather*}

Next, we estimate the dimensions of 
the indecomposable summands of $V^{\otimes n}_{M}$ which are from $ST_{r}(\Gamma_{M})$. We start with a general and well-known lemma.

\begin{Lemma}\label{l:glm}
We have $\dim T(a_{1},\dots,a_{M-1})\leq\prod_{i=1}^{M-1}\binom{M}{i}^{a_{i}}$.
\end{Lemma}

\begin{proof}
Recall that $\bigwedge^{i}V_{M}$ is a tilting $\Gamma_{M}$-representation for all $i\in\{1,\dots,M-1\}$
(this follows since $\bigwedge^{i}V_{M}$ is the Weyl representation 
$\Delta(0,\dots,0,1,0,\dots,0)$ for 
the $i$th fundamental weight $\omega_{i}$ and this weight is minimal in the set of dominant integral weights).
Now, essentially by their construction, the $\Gamma_{M}$-representation 
$T(a_{1},\dots,a_{M-1})$ is a direct summand of the $\Gamma_{M}$-representation
$(\bigwedge^{1}V_{M})^{\otimes a_{1}}\otimes\dots\otimes(\bigwedge^{M-1}V_{M})^{\otimes a_{M-1}}$, and $\dim\bigwedge^{i}V_{M}=\binom{M}{i}$.
\end{proof}

\begin{Lemma}\label{l:new}
Let $D_{n}$ denote the maximum of the dimensions of the indecomposable summands 
of $V^{\otimes n}_{M}$ from $ST_{r(n)}(\Gamma_{M})$. Then
\begin{gather*}
\lim_{n\to\infty}\sqrt[n]{D_{n}}=1.
\end{gather*} 
\end{Lemma}

\begin{proof}
Every such summand is of the form $St_{r}\otimes T^{(r)}$ where $T$ is an indecomposable tilting representation and $(\placeholder)^{(r)}$ is the $r$th Frobenius twist,
see \cite[Remark 2(1)]{An-steinberg-linkage}. Hence, the highest weight of $T$ should be bounded by $n/p^{r(n)}$ in the sense that sum of coefficients
of the fundamental weights is bounded by this number (more restrictively even, if $\lambda$ is the highest weight of $T$, then the weight $(p^r-1)\rho+p^r\lambda$ should appear in $V^{\otimes n}_{M}$).

If we let $A$ denote the maximum $A=\max_{i}\{\binom{M}{i}\}$, so that $A=\binom{M}{(M-1)/2}$, then by \autoref{l:glm}, we know that if the relevant tilting module in $ST_{r(n)}(\Gamma_{M})$ is to appear in $V^{\otimes n}_{M}$, then
\begin{gather*}
\dim St_{r}\otimes T^{(r)}\leq p^{r}A^{\sum_i a_i}\leq p^{r}A^{n/p^{r}}.
\end{gather*}
Now we can calculate
\begin{gather*}
\lim_{n\to\infty}\sqrt[n]{p^{r(n)} A^{n/p^{r(n)}}}\leq\lim_{n\to\infty}\sqrt[n]{\sqrt{n} A^{\sqrt{n}}}=1,
\end{gather*}
which concludes the proof.
\end{proof}

Now the total number of summands 
of $V^{\otimes n}_{M}$ coming from $ST_{r}(\Gamma_{M})$ is at least
$\frac{t_{n}}{D_{n}}$. Hence, 
\begin{gather*}
t_{n}/D_{n}\leq b_{n}^{M}\leq M^{n}
,\quad
\lim_{n\to\infty}\sqrt[n]{t_{n}/D_{n}}=M.
\end{gather*}
Hence, \autoref{p:glm} follows for $M$ odd. Then \autoref{l:avoideven} implies \autoref{p:glm} for $M$ even.

%%%%%%%%%%%%%%%%%%%%%%%%%%%%%%%%%%%%%%%%%%%%%%

\section{The general linear super group and consequences}\label{s:super}

%%%%%%%%%%%%%%%%%%%%%%%%%%%%%%%%%%%%%%%%%%%%%%

In this section we will work 
in the category of super vector spaces over $\bk$ (although we 
sometimes omit the word `super'` 
to avoid too cumbersome phrasings). Since the latter reduces to the ordinary category of vector spaces in characteristic 2, we assume $\mathrm{char}(\bk)\neq 2$.

Recall from \autoref{d:schemes} that
an affine group superscheme (in short: a supergroup) $G$ over $\bk$ is a representable functor from the category of commutative superalgebras (associative $\mZ/2\mZ=\{\oa,\ob\}$-graded algebras which are graded commutative) over $\bk$ to the category of groups. For general background on the theory of supergroups we refer to, for example, \cite{BrKlQ}, \cite{BrKu-mullineux-conjecture}, \cite{Masuoka} and \cite{Musson}.

We refer to $(M,N)\in\mN^{\times 2}$ as the `super dimension' of the $\mZ/2\mZ$-graded vector space $\bk^{M|N}$, and $M+N$ as the `dimension'. This should not lead to confusion as we will have no need for the `categorical dimension' $M-N$, also something referred to as the (super) dimension.

\begin{Proposition}\label{ThmSuper}
For a representation $V$ of a supergroup $G$ on a super vector space $V$ of super dimension $(M,N)$, we have
\begin{gather*}
\beta^{G,V}=\lim_{n\to\infty}\sqrt[n]{b_{n}^{G,V}}=M+N.
\end{gather*}
\end{Proposition}

Here the numbers $\beta^{G,V}$ and $b_{n}^{G,V}$ have the same meaning as before, i.e. they refer to the number of indecomposable summands in the $G$-representation $V^{\otimes n}$.
If we denote the representing commutative Hopf superalgebra for~$G$ by~$\cO(G)$, then a representation of $G$ can either be interpreted as a $\mZ/2\mZ$-graded comodule for $\cO(G)$, or equivalently as a homomorphism $G\to GL_{M|N}$ of supergroups.
By the latter interpretation, it is clearly sufficient to prove \autoref{ThmSuper} for $G=GL_{M|N}$ and $V=V_{M|N}$ its vector representation on $\bk^{M|N}$.

Before getting to the proof of \autoref{ThmSuper}, we derive some consequences.

\begin{proof}[Proof of \autoref{t3:main}.(b)]
If $\mathrm{char}(\bk)=0$, then
\begin{gather*}
\bk S_{n}\to\End_{GL_{M|N}}(V_{M|N}^{\otimes n})
\end{gather*}
is surjective, see \cite[Section 4]{BeRe}. We can therefore repeat the proof of \autoref{t3:main}.(a) from \autoref{s:groups} to reduce to \autoref{ThmSuper} for $G=GL_{M|N}$. If one wants to write things out explicitly, the set of partitions is now those for which $\lambda_{M+1}\leq N$.
\end{proof}

\begin{Remark}\label{RemSupP}
The proof of \autoref{t3:main}.(b) does not extend to positive characteristic. Indeed, in this case $\bk S_{n}\to\End_{GL_{M|N}}(V_{M|N}^{\otimes n})$ need not be surjective, see \cite[Theorem C]{CEKO}. More concretely, it is observed in \cite[\S 4]{CEKO}, that for $p=3$, $M=2$ and $N=1$, the number of indecomposable summands in $V^{\otimes 5}$ is 17, while the number of primitive idempotents in a decomposition of unity in $\bk S_{5}$ is only 16. Hence the action of the symmetric group on tensor powers of the vector representation of $GL_{M|N}$ is not sufficient to account for all indecomposable summands.
\end{Remark}

Because of this remark, we need an alternative proof for semigroups compared to the non-super case:

\begin{proof}[Proof of \autoref{t2:main}.(a) -- affine semigroup superschemes]
A representation of an affine semigroup superscheme $\Gamma$ corresponds to a semigroup homomorphism $\Gamma\to Mat_{M|N}$, with $Mat_{M|N}$ denoting the monoid superscheme of square $(M+N)$-matrices. In particular, the number of summands in $V^{\otimes n}$ over $\Gamma$ is bounded from below by the number of summands over $Mat_{M|N}$. By considering $\cO(Mat_{M|N})$ as a subcoalgebra of $\cO(GL_{M|N})$, we can identify the category of $Mat_{M|N}$-representations with the category of polynomial $GL_{M|N}$-representations, so the number of` direct summands in $V^{\otimes n}$ over $Mat_{M|N}$ is the same as over $GL_{M|N}$.

In conclusion, the number of direct summands over $\Gamma$ is bounded from below by the number of summands over $GL_{M|N}$. The result thus follows from \autoref{ThmSuper} for $G=GL_{M|N}$.
\end{proof}

\subsection{Proof of \autoref{ThmSuper} -- semisimple case}

Assume that $\mathrm{char}(\bk)=0$. We get a stronger statement:

\begin{Lemma}\label{l:supersemisimple}
The $GL_{M|N}$-representation $V_{M|N}^{\otimes n}$ is semisimple and the dimension of the simple representations occurring in $V_{M|N}^{\otimes n}$ is bounded by a polynomial in $n$ of degree $\frac{M(M-1)+N(N-1)}{2}$.
\end{Lemma}

\begin{proof}
That the tensor powers are semisimple is proved in \cite[Theorem 5.14]{BeRe}. The dimension of these simple representations is bounded by that of the Kac modules with same highest weight, see for instance \cite[\S 8.2]{Musson}. As induced modules, the dimension of the latter is given by a constant (depending on $MN$, see for instance \autoref{LemInd}(b)) times the dimension of the simple $(GL_{M}\times GL_{N})$-representation with same highest weight, which is a weight appearing in $V_{M|N}^{\otimes n}$. The latter can be bounded by a polynomial in $n$ of degree $\frac{M(M-1)+N(N-1)}{2}$, as explained in \autoref{ss:semisimple}.
\end{proof}

Let $b_{n}^{M,N}$ be the analog of $b_{n}^{M}$ for $GL_{M|N}$.

\begin{Proposition}\label{p:supersemisimple}
We have
\begin{gather*}
\Omega\big((M+N)^{n}/n^{(M(M-1)+N(N-1))/2}\big)\ni b_{n}^{M,N}\in O^{\prime}\big((M+N)^{n}\big).
\end{gather*}
\end{Proposition}

\begin{proof}
As before, this follows from \autoref{l:supersemisimple} and the super 
analog of \autoref{l:main}.
\end{proof}

\subsection{Preparation for the proof: distributions and induction}\label{SecPrep}

By a `subgroup' $H<G$ of a supergroup we refer to a representable subgroup functor, or equivalently a closed subsuperscheme which is also closed under the group operation.

{For the general linear supergroup $GL_{M|N}$ we consider the subgroups $P^{+},P^-< GL_{M|N}$.} Here, for any superalgebra $A$
\begin{gather*}
P^{+}(A)<GL_{M|N}(A)=\Aut_{A}(A^{M|N})
\end{gather*}
consists of all automorphisms which are expressed as ($M+N$)-block matrices in a way that the left down block of size $N\times M$ is zero. The subgroup $P^{-}$ corresponds similarly to a zero ($M\times N$)-block.

For a supergroup $G$, we have the underlying affine group scheme $G_{0}$, which can be defined as the restriction of the functor $G$ to $\bk$-algebras (viewed as superalgebras contained in degree $\oa$) or via the quotient of $\cO(G)$ by the ideal generated by all odd elements. For $G=GL_{M|N}$ we have $G_{0}=P^{+}\cap P^{-}=GL_M\times GL_{N}$.

For an affine group scheme $G$, one defines the distribution superalgebra  as a subalgebra $\Dist G\subset\cO(G)^{\ast}$ which is a cocommutative Hopf superalgebra, similarly to the classical case, see \cite[\S 3]{BrKu-mullineux-conjecture}. Explicit descriptions of these algebras for $GL_{M|N}$ and subgroups as $P^{\pm}$ are also given loc. cit.

We also have the Lie superalgebra $\Lie G$ as a subspace of $\Dist G$. For $GL_{M|N}$ this is the general linear Lie superalgebra $\mathfrak{gl}_{M|N}$ of square ($M+N$)-matrices with supercommutator. We have a vector space decomposition
\begin{gather*}
\mathfrak{gl}_{M|N}=\fg_{-}\oplus\fg_{\oa}\oplus\fg_{+},
\end{gather*}
where $\fg_{\oa}\oplus \fg_{\pm}$ is the subalgebra corresponding to $P^{\pm}< G$.

For a supergroup $G$ we denote by $\Rep_{\bk}G$ the rigid monoidal category of finite dimensional (super) representations. In particular, for $G$ an ordinary affine group scheme interpreted as a supergroup, this category is equivalent (as a $\bk$-linear additive category) to a direct sum of two copies of the classical representation category.

\begin{Lemma}\label{LemInd}
Set $G=GL_{M|N}$.
\begin{enumerate}

\item The forgetful functors $\res^{G}_{P}\colon\Rep_{\bk}G\to\Rep_{\bk}P$ and $\res^{P^{-}}_{G_{0}}\colon\Rep_{\bk}P^{-}\to\Rep_{\bk}G_{0}$ have left adjoint functors $\ind^{G}_{P}$ and $\ind^{P^{-}}_{G_{0}}$.

\item As $G_{0}$-representations and $\fg_{-}$-representations, we have
\begin{gather*}
\ind^{G}_{P}M\;\simeq\;\Lambda\fg_{-}\otimes M.
\end{gather*}

\item We have a natural isomorphism 
\begin{gather*}
\ind^{P^{-}}_{G_{0}}\res^{P}_{G_{0}}\;
\Rightarrow\;\res^{G}_{P^{-}}\ind^{G}_{P}.
\end{gather*}

\end{enumerate}
\end{Lemma}

\begin{proof}
This can be proved by relying either on the theory of Harish-Chandra pairs from \cite{Masuoka} or the distribution algebras from \cite{BrKu-mullineux-conjecture}. We choose the latter approach.

In \cite{BrKu-mullineux-conjecture} it is proved that, for $H$ denoting any of the supergroups in the lemma, $\Rep_{\bk}H$ is equivalent to the category of integrable finite dimensional modules of $\Dist H$. Here, `integrable' essentially means weight module. Clearly, on the level of ($\Dist H$)-representations, restriction has a left adjoint functor given by induction, for instance $\Dist G\otimes_{\Dist P}{}_{-}$. By the explicit realization in \cite[\S 4]{BrKu-mullineux-conjecture}, it follows that 
\begin{gather*}
\Dist G\;\simeq\;\Lambda\fg_{-}\otimes\Dist P,\quad\text{respectively}\quad\Dist P^{-}\;\simeq\;\Lambda\fg_{-}\otimes\Dist G_{0}
\end{gather*}
as right ($\Dist P$)-representations respectively right ($\Dist G_{0}$)-representations. It follows easily that induction sends integrable modules to integrable modules, providing the desired left adjoints in (a). Statements (b) and (c) then follow again from the above and the explicit forms of the distribution algebras in \cite[\S 4]{BrKu-mullineux-conjecture}.
\end{proof}

\begin{Lemma}\label{LemGH}
Consider a supergroup $G$ with subgroup $H<G$ and a representation $V$ of $G$ on $\bk^{M|N}$. Assume that:
\begin{enumerate}[label=\emph{\upshape(\roman*)}]

\item $\beta^{H,V}=M+N$.

\item $\res^{G}_H\colon\Rep_{\bk}G\to\Rep_{\bk}H$ has a left adjoint $\ind^{G}_{H}$ and
there exist $j\in\mN$ and $U\in\Rep_{\bk}H$ for which $\ind^{G}_{H}U$ is a $G$-summand of $V^{\otimes j}$.

\end{enumerate}
Then it follows that
$\beta^{G,V}=M+N$.
\end{Lemma}

\begin{proof}
As a direct consequence of adjunction and
the definition of $\ind^G_H$, we find for any $W\in\Rep_{\bk}G$
\begin{gather*}
\ind^G_H(U)\otimes W\;\simeq\;\ind^G_H\big(U\otimes\res^G_H(W)\big).
\end{gather*}
By assumption, $b_{n+j}^{G,V}$ is at least the number of $G$-summands in
\begin{gather*}
\ind^G_H (U)\otimes V^{\otimes n}\;\simeq\;\ind^G_H(U\otimes V^{\otimes n}),
\end{gather*}
which in particular shows that $b_{n}^{H,V}\leq b_{n+j}^{G,V}$.
Thus, assumption (i) and the super 
analog of \autoref{l:main} imply the claim.
\end{proof}

\subsection{Proof of \autoref{ThmSuper} -- positive characteristic}

Now we fix $M,N>0$, $G=GL_{M|N}$, $P^{\pm}<G$ from \autoref{SecPrep} and thus $G_{0}=P_{0}=GL_{M}\times GL_{N}$. Let $V=\bk^{M|N}$ be the vector representation of $G$. Since we will use comparison with characteristic zero, we do not yet make assumptions on $\mathrm{char}(\bk)$.

For partitions $\lambda=(\lambda_{i})_{1\leq i\leq M}$, $\mu=(\mu_{j})_{1\leq j\leq N}$ of length at most $M$, $N$, we denote by $L_{0}(\lambda|\mu)$ the corresponding simple polynomial $G_{0}$-representation, which we also interpret as a $P$-representation in the usual way (for instance with trivial action of $\fg_+$).

\begin{Lemma}\label{LemTyp}
The $G$-representation
$\ind^G_{P} L_{0}(\lambda|\mu)$
is simple if the integer
\begin{gather*}
c_{ij}(\lambda|\mu):=\lambda_{i}-i+\mu_{j}-j+M+1
\end{gather*}
is not zero in $\bk$, for all $1\leq i\leq M$ and $1\leq j\leq N$.
\end{Lemma}

\begin{proof}
By \autoref{LemInd}(b), $\ind^G_P L_{0}(\lambda|\mu)$ is free as a $\Lambda\fg_{-}$-representation, so every $P^{-}$-submodule contains $u\otimes L_{0}(\lambda|\mu)$ for $u$ a non-zero element (unique up to constant) in the top degree of $\Lambda\fg_{-}$. Let $v_{+}$ be a highest weight vector of $L_{0}(\lambda|\mu)$. It suffices to prove that the $\fg_{+}$-submodule of $\ind^G_P L_{0}(\lambda|\mu)$ generated by $u\otimes v_{+}$ contains $1\otimes v_{+}$, where we use $\ind^G_P L_{0}(\lambda|\mu)\simeq\Lambda\fg_{-}\otimes L_{0}(\lambda|\mu)$. By choosing a conveniently ordered product of all root vectors in $\fg_{-}$ for $u$ and a mirrored product of root vectors for $\fg_{+}$ for $v$, it follows that
\begin{gather*}
vu\otimes v_{+}=\prod_{ij}c_{ij}(\lambda|\mu) \;(1\otimes v_{+}),
\end{gather*}
which concludes the proof.
\end{proof}

Now we fix a prime $p$ and consider the partitions $\alpha$, $\nu$ of lengths $M$, $N-1$ given by\begin{gather*}
\alpha_{i}
=N+(p-1)(M-i),1\leq i\leq M,
\quad
\nu_{j}
=(p-1)(N-j),1\leq j\leq N.
\end{gather*}
Since $\alpha_M=N$ is greater than the length of $\nu$, the partition $\kappa:=\alpha\nu^{t}$ of length $M+\nu_{1}$ makes sense. Concretely
\begin{gather*}
\kappa_{i}=
\begin{cases}
\alpha_{i} & \mbox{ if }i\leq M
\\
(\nu^{t})_{i-M} & \mbox{ if }i>M.
\end{cases}
\end{gather*}
It follows that $\kappa$ is a $p$-core. 

Before coming to a crucial 
proposition, we need the following (well-known) lemma.

\begin{Lemma}\label{L:core}
Let $\lambda$ be a $p$-core partition, and $T$ be a standard $\lambda$-tableaux. The primitive idempotent $e_{T}$ in $\BQ S_{r}$ associated to $T$ via Young symmetrizers (so that in particular $(\BQ S_{r})e_{T}$ is isomorphic to the simple Specht module $S^\lambda$) belongs to $\mZ_{(p)}S_{r}\subset\BQ S_{r}$.
\end{Lemma}

\begin{proof}
According to the theory of Young symmetrizers, see for example \cite[Section 7.2]{Fu-tableaux},
we can clear denominators 
in $e_{T}$ and get a pseudo-idempotent $\tilde{e}_{T}\in\BZ S_{r}$. This pseudo-idempotent satisfies $\tilde{e}_{T}\tilde{e}_{T}=n_{\lambda}\cdot\tilde{e}_{T}$ (in other words $e_T=\tilde{e}_T/n_{\lambda}$), where $n_{\lambda}\in\mZ$ is the product of the {hook lengths} in $\lambda$ by e.g. \cite[Section 7.4, Exercises 18 and 19]{Fu-tableaux}. Moreover, the {hook lengths} of $p$-cores are never divisible by $p$, see e.g. 
\cite[Statement 2.7.40]{JaKe-reptheory-symmetric-group}.
Hence, for a $p$-core we have $n_{\lambda}\notin p\BZ$, and we are done.
\end{proof}

For the remainder of the section, we assume that $\mathrm{char}(\bk)=p>2$.
By the above lemma, we can take an idempotent in $e_0\in\BZ_{(p)}S_r$, with $r=|\kappa|$, which is a primitive idempotent corresponding to $\kappa$, when considered in $\BQ S_r$, and we denote by $e\in \BF_p S_r\subset \bk S_r$ its image modulo $p$. We denote by $\mS^{\kappa}V$ the corresponding direct summand $e(V^{\otimes r})$ in $V^{\otimes r}$ and we use the same notation $\mS^{\kappa}V$ for the summand $e_0(V^{\otimes r})$ when working over $\BQ$.

\begin{Proposition}\label{MainProp}
Working either over $\bk^{\prime}=\BQ$ or $\bk^{\prime}=\bk$ and with $\alpha$, $\nu$, $\kappa$ as just introduced, we have
\begin{gather*}
\mS^{\kappa}V\;\simeq\;\ind^{G}_{P}L_{0}(\alpha|\nu)\qquad\mbox{in }\; \,\Rep_{\bk^{\prime}}GL_{M|N}.
\end{gather*}
\end{Proposition}

\begin{proof}
By construction (and realizing $V$ as the extension of scalars over a free $\mZ_{(p)}$-module), the character of the left-hand side is identical for $\bk$ and $\BQ$. The $GL_{M}\times GL_{N}$-representation $L_{0}(\alpha|\nu)$ has the Weyl character, because $\nu$ is a $p$-core and $\alpha$ is just a `shifted' $p$-core. Indeed, this follows from \cite[Theorem~4.5]{JaMa-q-jantzen-schaper}, or from a direct application of tilting theory and \autoref{L:core}.
It therefore follows from \autoref{LemInd}(2) that also the character of $\ind^{G}_{P}L_{0}(\alpha|\nu)$ is identical over $\bk$ and $\BQ$.

By \autoref{LemTyp}, the right-hand side is a simple representation. Indeed, we can calculate
\begin{gather*}
c_{ij}(\alpha|\nu)=(M+N-i-j)p+1,
\end{gather*}
which is never zero in $\bk$ or in $\BQ$. 
In particular, over $\BQ$, the right-hand side is the simple representation with highest weight $\alpha|\nu$. The isomorphism over $\BQ$ now follows from \cite[Section 4]{BeRe}.

Now working over $\bk$, by the combination of the previous two paragraphs, the two representations have the same character and the right-hand side is simple. They must thus be isomorphic.
\end{proof}

Now we can prove the main result.

\begin{proof}[Proof of \autoref{ThmSuper}]
By the combination of \autoref{MainProp} and \autoref{LemGH} it suffices to prove \autoref{ThmSuper} for the supergroup $P$ acting on $V$. We prove the equivalent formulation in terms of $P^{-}$.

\autoref{MainProp} and \autoref{LemInd}(c) imply that
\begin{gather*}
\mS^{\kappa}V\;\simeq\;\ind^{P^{-}}_{G_{0}}L_{0}(\alpha|\nu)
\end{gather*}
as $P^{-}$-representations. In particular, the result for $(P^{-},V)$ follows, via \autoref{LemGH} from the purely even case $(G_{0},V)$ in \autoref{s:groups}.
\end{proof}

%%%%%%%%%%%%%%%%%%%%%%%%%%%%%%%%%%%%%%%%%%%%%%

\section{The general linear quantum group and consequences}\label{s:qgroups}

%%%%%%%%%%%%%%%%%%%%%%%%%%%%%%%%%%%%%%%%%%%%%%

Let $\bk$ still denote an arbitrary field. Further, 
fix $q\in\bk^{\ast}$. 
We consider either of the following objects.

To $(\bk,q)$ we associate the pair, often called the mixed characteristic of $(\bk,q)$, by
\begin{gather*}
(p,\ell):=(|1|,|q^{2}|)\in 
(\BN\cup\{\infty\})^{\times 2}
\end{gather*}
of the orders $|\placeholder|$ of $1$ and $q^{2}$ in the additive group underlying $\bk$. For example, the case $p=\ell\in\mN$ corresponds to a field $\bk$ of positive characteristic and $q=1$ and the case $p=\ell=\infty$ corresponds to a field of characteristic zero with generic $q$.
We also use the quantum numbers for $a\in\BN$ and $x\in\bk$:
\begin{gather*}
[a]_{x}=x^{-(a-1)}+x^{-(a-3)}+\dots+x^{a-3}+x^{a-1}=\frac{x^{a}-x^{-a}}{x-x^{-1}}\in\bk,
\end{gather*}
where the second equality is only applicable for $q^{2}\neq1$. It then follows that $\ell$ is also the minimal value of $n$ for which $[n]=0$, or $\infty$ when no such value exists.

We call $q\in\bk^{\ast}$ generic if $q$ is not a root of unity in 
$\overline{\bk}$, e.g. $q$ could be the formal variable in the field $\BC(q)$ of rational complex functions.

\begin{Example}\label{e:mixedchar}
Let us give a few examples:
\begin{enumerate}

\item For $\bk=\BC(q)$ and generic $q$ we have $(p,\ell)=(\infty,\infty)$.
For $\bk=\BF_{7}(q)$ and generic $q$ we have $(p,\ell)=(7,\infty)$.
These cases are both `semisimple', in the terminology of \autoref{RemMix}.

\item For $\bk=\BC$ and $q=\exp(\pi i/3)$ we have $(p,\ell)=(\infty,3)$.

\item For $\bk=\BF_{7}$ and $q=2$ we have $(p,\ell)=(7,3)$.

\item For $\bk=\BF_{7}$ and $q=1$ we have $(p,\ell)=(7,7)$.

\end{enumerate}
In general, 
the representation category of $\Gamma_{M}$ is semisimple if and only if 
$\ell=\infty$.
\end{Example}

We consider Lusztig's divided power quantum group over $\bk$ as in \cite{Lu-qgroups-root-of-1} associated to a type $A$ Cartan datum. 
We use $\Gamma_{M}=\Gamma_{M}=U_{q}^{\bk}(\mathfrak{sl}_{M})$ or 
$\Gamma_{M}=U_{q}^{\bk}(\mathfrak{gl}_{M})$ as the notation.

\begin{Remark}\label{RemMix}
Following \autoref{e:mixedchar}, all of our discussions regarding $\Gamma_{M}$ split into 
four distinct cases:
\begin{enumerate}

\item For $(p,\infty)$ the quantum group representations are semisimple for any $p$, and their combinatorics is the same as for $G(\BC)$, for $G=SL_M$ or $G=GL_M$.
We call this the semisimple case. We discuss this case in \autoref{ss:qsemisimple}.

\item The case $(\infty,\ell)$ for $\ell<\infty$
can be combinatorially identified with its special case $\bk=\BC$. 
We call this the complex 
quantum group case. We discuss this case in \autoref{ss:qproofglcomplex}.

\item The strictly mixed case is $p,\ell<\infty$ and $p\neq\ell$. We discuss this case in 
\autoref{ss:qproofglm}.

\item The 
situation $p=\ell<\infty$ prime is characteristic $p$. This reduces to the case in \autoref{s:groups}.

\end{enumerate}
This list is ordered in increasing order of difficulty, in the sense that
the dimensions of the indecomposable summands are increasing, reading 
from top to bottom (and also more difficult to compute). Hence, using the same definitions as before, the convergence rate of $\sqrt[n]{b_{n}^{\Gamma_{M},V}}$ is slower 
for the bottom cases compared to the top ones.
\end{Remark}

Let $V_{M}$ denote the quantum vector representation of 
$\Gamma_{M}$. With the same notation as in the previous sections we have the analog of \autoref{p:glm}:

\begin{Proposition}\label{p:glmq}
For any $M\in\BZ_{>0}$ we have 
\begin{gather*}
\lim_{n\to\infty}\sqrt[n]{b_{n}^{M}}=M.
\end{gather*}
\end{Proposition}

\begin{proof}
The case $M=1$ is again immediate, and $M>1$ is proven
in \autoref{ss:qproofgltwo} and \autoref{ss:qproofglm}.
\end{proof}

\begin{Remark}\label{r:qtilting}
As before, we make use of the fact that $V^{\otimes n}_{M}$ is tilting, see e.g. \cite[Proposition 2.3]{AnStTu-semisimple-tilting}. In particular, $V^{\otimes n}_{M}$ is a direct sum 
of indecomposable tilting representations.
For $(m_{1},\dots,m_{M-1})\in\BN^{M-1}$ we use the Weyl representations 
$\Delta(m_{1},\dots,m_{M-1})$
and the indecomposable tilting representations $T(m_{1},\dots,m_{M-1})$ 
of highest weight $(m_{1},\dots,m_{M-1})$.
\end{Remark}

We now consider any $\bk$-subalgebra $\Gamma\subset U_{q}^{\bk}(\mathfrak{gl}_{M})$. In this case $V^{\otimes n}$ is a $\Gamma$-representation by restriction, although the tensor product of 
representations does not need to exist in general.

\begin{proof}[Proof of \autoref{t2:main} -- quantum groups of type $A$ and $\bk$-subalgebras]
The only difference to the proof in \autoref{s:groups} 
is that we use \autoref{p:glmq} instead of \autoref{p:glm}.
\end{proof}

\begin{Remark}\label{r:coideal}
As observed in the 1990s or even earlier, 
embeddings of Lie subalgebras $\mathfrak{g}\hookrightarrow\mathfrak{gl}_{M}$ 
do not quantize properly. Consequently, contrary to \autoref{s:groups} and~\autoref{s:super}, proving \autoref{t2:main} for the general linear case does not imply it for other quantum groups.
For this reason, quantum groups that are not of type 
$A$ are outside of the scope of this paper.
On the other hand, \autoref{t2:main} does include 
coideal subalgebras, with the 
most prominent example being quantum symmetric pairs (also called 
$\imath$quantum groups), which have been studied many people, see e.g. \cite{NoSu-q-symmetric-spaces}, 
\cite{Le-symmetric-pairs} or \cite{Ko-qsym-kac-moody}.
\end{Remark}

%%%%%%%%%%%%%%%%%%%%%%%%%%%%%%%%%%%%%%%%%%%%%%

\subsection{Proof of \autoref{p:glmq} -- semisimple case}\label{ss:qsemisimple}

%%%%%%%%%%%%%%%%%%%%%%%%%%%%%%%%%%%%%%%%%%%%%%

In this case, the indecomposable tilting representation 
$T_{q}(m_{1},\dots,m_{M-1})$ is isomorphic to
the Weyl representation $\Delta_{q}(m_{1},\dots,m_{M-1})$, and the latter 
has the quantum Weyl character, see e.g. \cite[Equation 2]{Sa-quantum-roots-of-unity}.  We therefore get the same bound 
as for the group schemes, namely $\frac{M^{n}}{f(n)}\leq b_{n}^{M}$ for
$f\in\BN[x]$ with $\deg f=\tfrac{M(M-1)}{2}$.
As before, an analog 
of \autoref{p:semisimple} follows from the above. 
This in turn implies \autoref{p:glmq} and \autoref{t2:main}. 

\begin{Example}\label{e:qclt}
Let $\Gamma_{2}=U_{q}^{\bk}(\mathfrak{sl}_{2})$ for $(p,\infty)$.
As in \autoref{e:clt}, we see that 
$b_{n}^{G,V_{2}}\sim 2^{n}/(\pi n/2)^{1/2}$ where 
$V_{2}$ is the quantum vector representation.
Moreover, for any nontrivial representation of $\Gamma_{2}$ 
we get that
$b_{n}^{G,V}\sim(\dim V)^{n}/(2\pi nA^{2})^{1/2}$ 
for some $A=A(V)\in\BR_{>0}$.
\end{Example}

%%%%%%%%%%%%%%%%%%%%%%%%%%%%%%%%%%%%%%%%%%%%%%

\subsection{Proof of \autoref{p:glmq} -- for \texorpdfstring{$M=2$}{M=2}}\label{ss:qproofgltwo}

%%%%%%%%%%%%%%%%%%%%%%%%%%%%%%%%%%%%%%%%%%%%%%

First, in the case $M=2$ we can use the known character formulas 
for $\Gamma_{2}$:

\begin{Example}\label{e:qsl2}
The quantum version of \autoref{e:sl2} is as follows, 
see e.g. \cite[Section 2]{SuTuWeZh-mixed-tilting} for details.
We will use the above quantum numbers for $x=v$ 
a formal variable for the characters.

Let $p^{(i)}=p^{i-1}\ell$ for $i>0$ and $p^{(0)}=1$.
Let $m+1=a_np^{(n)}+\cdots +a_{1}p^{(1)}+a_{0}p^{(0)}=(a_{n},\dots,a_{0})$
the $(p,\ell)$-adic expansion of $m+1$ where $a_{0}\in\{0,1,\dots,\ell-1\}$ and $a_{i}\in\{0,1,\dots,p-1\}$ 
for $i>0$. 
(Note that the cases $p=\infty$ or $\ell=\infty$ are covered by the notation as well.)
We then have
\begin{gather*}
\mathrm{ch}\,T(m)=
[a_{n}p^{(n)}]_{v}\cdot\prod_{a_{i}\neq 0,i\neq n}
[2]_{v^{a_{i}p^{(i)}}}.
\end{gather*}
For example, for $m=52$ and $p=2$ and $\ell=3$ we have 
$m+1=(1,0,0,0,1,2)$ so
\begin{gather*}
\mathrm{ch}\,T(52)=
[p^{(5)}]_{v}\cdot 
[2]_{v^{p^{(1)}}}[2]_{v^{2p^{(0)}}}.
\end{gather*}
For $v=1$ we get $\dim T(52)=192$.
\end{Example}

This proves:

\begin{Lemma}\label{l:qsl2}
Let $\alpha=1$ for $\ell=\infty$, and otherwise let $\alpha=1+(\log_{2} p^{\prime})^{-1}$ 
where $p^{\prime}=\min\{p,\ell\}$. Then for any $m\leq n$ we have
\begin{gather*}
\dim T(m)\leq(n+1)^{\alpha}.
\end{gather*}
\end{Lemma}

\begin{proof}
By \autoref{e:qsl2}, as in the proof of \autoref{l:sl2}.
\end{proof}

As before, we get the analog of \autoref{p:sltwo} from \autoref{l:qsl2}.
This in turn implies \autoref{p:glmq}.

\begin{Example}
Let $m=52$, so $m+1=53$. Here are some $(p,\ell)$-adic expansions 
and the dimension of $T(m)$:
\begin{gather*}
(p,\infty)\colon 53=(53),\quad\dim T(52)=53,
\\
(\infty,3)\colon 53=(17,2),\quad\dim T(52)=102=2^1\cdot 17\cdot 3,
\\
(2,3)\colon 53=(1,0,0,0,1,2),\quad\dim T(52)=192=2^2\cdot(2^4\cdot 3),
\\
(2,2)\colon 53=(1,1,0,1,0,1),\quad\dim T(52)=256=2^3\cdot 2^5.
\end{gather*}
Note that the smaller $(p,\ell)$ the bigger $\dim T(m)$.
\end{Example}

With $\alpha$ as in \autoref{l:qsl2} we then again get 
$\frac{2^{n}}{(n+1)^{2}}\leq\frac{2^{n}}{(n+1)^{\alpha}}\leq b_{n}^{2}\leq 2^{n}$, 
which again proves \autoref{t2:main} and gives a bit finer information.

%%%%%%%%%%%%%%%%%%%%%%%%%%%%%%%%%%%%%%%%%%%%%%

\subsection{Proof of \autoref{p:glmq} -- for \texorpdfstring{$M\geq 2$}{M>=2}}\label{ss:qproofglm}

%%%%%%%%%%%%%%%%%%%%%%%%%%%%%%%%%%%%%%%%%%%%%%

We now show that \autoref{p:glmq} holds for $(p,\ell)$ with $p,\ell<\infty$
and $p\neq\ell$, and arbitrary $M\geq 2$.

We can use the same arguments as in the reductive group case 
above, with the following adaptations (the reader should compare \autoref{e:sl2} and 
\autoref{e:qsl2} while reading the below):

\begin{enumerate}

\item Recall that $p^{(i)}=p^{i-1}\ell$ for $i>0$ and $p^{(0)}=1$. All appearances of $p^{r}$ should be replaced by $p^{(r)}$.

\item The results we need from \cite{An-steinberg-linkage} hold, mutatis mutandis, 
for the quantum group as well, see \cite[Remark 2.2]{An-steinberg-linkage}.
That is, instead of the Frobenius twist 
one uses the Frobenius--Lusztig twist which, roughly speaking, acts as the Frobenius twist on digits $a_{i}$ for $i>0$ and as its quantum analog on the zeroth digit, and the rest is the same.

\end{enumerate}

Taking all of the above together, \autoref{p:glmq} follows 
from the same arguments as for $SL_{M}$, which proves
\autoref{t2:main}.

%%%%%%%%%%%%%%%%%%%%%%%%%%%%%%%%%%%%%%%%%%%%%%

\subsection{Some extra observations for the complex quantum group case}\label{ss:qproofglcomplex}

%%%%%%%%%%%%%%%%%%%%%%%%%%%%%%%%%%%%%%%%%%%%%%

For $p=\infty$ the combinatorics (in particular, the multiplicities 
of the decompositions) are the same as for the complex root of unity case.
Here, we consider $(p=\infty,\ell<\infty)$, since the case 
$(\infty,\infty)$ is semisimple and has the same combinatorics as the complex group case. This case is therefore already addressed in 
\autoref{ss:semisimple} and \autoref{ss:qsemisimple}.

\begin{Proposition}\label{p:qslm}
For $(m_{1},\dots,m_{M-1})\in\BN^{M-1}$.
Then there exists $A\in\BR_{>0}$ such that
\begin{gather*}
\dim T(m_{1},\dots,m_{M-1})\leq A\cdot\dim\Delta(m_{1},\dots,m_{M-1})
\end{gather*}
which implies
\begin{gather*}
\frac{M^{n}}{f(n)}\leq b_{n}^{M},
\quad f\in\BN[x],\deg f=\tfrac{M(M-1)}{2}
\end{gather*}
and
\begin{gather*}
\Omega\big(M^{n}/n^{M(M-1)}\big)\ni b_{n}^{M}\in O^{\prime}(M^{n}).
\end{gather*}
\end{Proposition}

\begin{proof}
In this case the numbers (tilting:Weyl)
are known to be given by parabolic Kazhdan--Lusztig polynomials, see \cite{So-tilting-a} and \cite{So-tilting-b}.
Even better, for $\Gamma_{M}=SL_{M}$ the results of 
\cite{St-diplom} imply that the parabolic Kazhdan--Lusztig polynomials are bounded. This in turn implies, again using the quantum Weyl character formula 
as in \autoref{ss:qsemisimple}, that the dimension of the tilting representation $T(m_{1},\dots,m_{M-1})$ is a polynomial in $m_{1},\dots,m_{M-1}$ of degree $\frac{M(M-1)}{2}$, as before. 
\end{proof}

Thus, \autoref{p:qslm} implies that
\autoref{p:glmq} holds for $(\infty,\ell)$ and arbitrary $M$, but the result 
is even a bit stronger.
We can even say a little more for $M\in\{2,3\}$:

\begin{Proposition}\label{L:q2and3}
\leavevmode

\begin{enumerate}

\item Let $M=2$. 	
For $m\in\BN$ we have
\begin{gather*}
\dim T(m)\leq 2(m+1),\quad
\frac{M^{n}}{2(n+1)}\leq b_{n}^{M}.
\end{gather*}

\item Let $M=3$. For $(m_{1},m_{2})\in\BN^{2}$ we have
\begin{gather*}
\dim T(m_{1},m_{2})\leq 12(m_{1}+1)(m_{2}+1)(m_{1}+m_{2}+2)
,\quad
\frac{M^{n}}{2(n^{3}+3n^{2}+3n+1)}\leq b_{n}^{M}.
\end{gather*}

\end{enumerate}
\end{Proposition}

\begin{proof}
(a). By \autoref{e:qsl2}.

(b). We will crucially use that the tilting characters (the Weyl multiplicities with the indecomposable tiling representations) are known 
explicitly by, for example, \cite{So-tilting-a}, \cite{So-tilting-b} and
\cite{St-diplom}. This explicit description of the characters is known as periodic patterns.

The quantum version of Weyl character formula gives $\dim\Delta(m_{1},m_{2})=(m_{1}+1)(m_{2}+1)(m_{1}+m_{2}+2)$.
We will use this as follows. 
The periodic pattern for $U_{q}^{\bk}(\mathfrak{sl}_{3})$ tilting representations are given by
\begin{gather*}
\includegraphics[height=3.5cm]{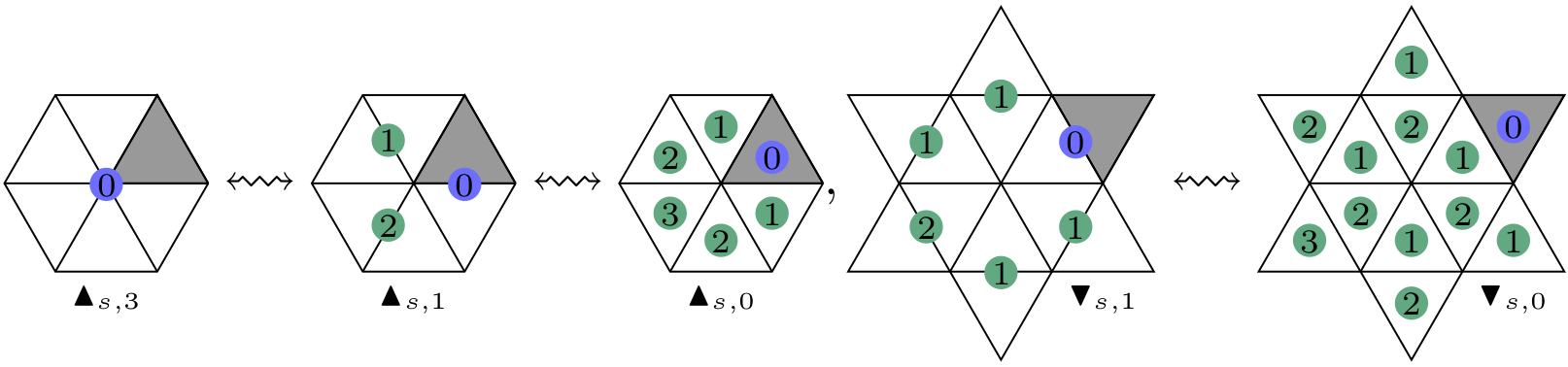}
\end{gather*}
These patterns mean that e.g. the tilting representation with highest
weight in a star pattern at the north east (the position of the highest weight is indicated by zero) has the twelve Weyl representations 
indicated by the circles in 
its Weyl filtration. All appearing highest weights 
of the Weyl representation are in the action orbit 
of the affine Weyl group on this alcove picture. Thus, 
$\dim\Delta(m_{1},m_{2})=(m_{1}+1)(m_{2}+1)(m_{1}+m_{2}+2)$ implies that 
$\dim T(m_{1},m_{2})\leq 12(m_{1}+1)(m_{2}+1)(m_{1}+m_{2}+2)$ for the star 
pattern. All other patterns have fewer Weyl factors and the claim follows.
The same bound, as one easily checks, works for the periodic patterns along the boundary of the Weyl alcove as well.
\end{proof}

Similar to \autoref{ss:qproofglm} we have certain summands that appear 
often enough to imply \autoref{t:main}:

\begin{Proposition}\label{l:q2monoid}
Let $M=2$ and $\ell=3$. Then there exist a family of summands 
of $V_{2}^{\otimes n}$ whose number of appearance $t_{n}$ 
in $V_{2}^{\otimes n}$ satisfies
\begin{gather*}
2^{n}/n^{5/2}\leq t_{n}.
\end{gather*}
\end{Proposition}

\begin{proof}
Via quantum Schur--Weyl duality, this is \cite[Theorem 4E.2]{KhSiTu-monoidal-cryptography}.
\end{proof}

\begin{Remark}\label{r:q2monoid}
The restriction to $\ell=3$ in \autoref{l:q2monoid} 
is used as \cite{KhSiTu-monoidal-cryptography} study the monoid 
version of the Temperley--Lieb calculus. Similar (and sharper) results 
can be obtained for any $\ell\in\BN$ by using Schur--Weyl 
duality and \cite[Propositions 9.4 and 9.5]{Sp-modular-tl}.
\end{Remark}

%%%%%%%%%%%%%%%%%%%%%%%%%%%%%%%%%%%%%%%%%%%%%%

\section{Counterexamples}\label{s:counter}

%%%%%%%%%%%%%%%%%%%%%%%%%%%%%%%%%%%%%%%%%%%%%%

In this section we will show that \autoref{t:main} does not extend arbitrarily, even over $\bk=\BC$.
We will use analog notion as before.

\begin{Theorem}\label{t:mainfalse}
Let $\bk=\BC$.
For every $A\in\BR_{\geq 0}$ there exists a cotriangular Hopf algebra $\Gamma$ with a 
finite dimensional $\Gamma$-corepresentation $V$ such that
\begin{gather*}
\beta^{\Gamma,V}+A<\dim V.
\end{gather*}
\end{Theorem}

\begin{proof}
The main player in this proof is the Temperley--Lieb category 
$TL(-2)$ of Rumer--Teller--Weyl \cite{RuTeWe-sl2}
with circle parameter $-2\in\BC$. It is the diagrammatic 
incarnation of $\cD=\Rep_{\mC}SL_{2}(\mC)$.

Let $m\geq 2$ and let $X=\BC^{2}$ be the 
vector representation of $SL_{2}(\mC)$. It follows from \cite[Theorem 1.1]{Bi-repcat-qgroup} that every matrix $E\in GL_{m}(\mC)$ with
\begin{gather*}
\mathrm{tr}(E^{T}E^{-1})=-2\in\mC
\end{gather*}
gives a nonsymmetric fiber functor
\begin{gather*}
F_{E}\colon\cD\to\cVect[\mC]
,
\quad X\mapsto\BC^{m}.
\end{gather*}
Here we use the notion fiber functor as in e.g. \cite[Definition 5.1.1]{EtGeNiOs-tensor-categories}.

Since the number of indecomposable summands 
$b_{n}^{\cD,X}$ does not depend on 
$F_{E}$ but only on $\Rep_{\mC}SL_{2}(\mC)$, the above implies that
\begin{gather*}
\beta^{\cD,X}=2<\dim X=m.
\end{gather*}

Finally, \cite[Theorem 1.1]{Bi-repcat-qgroup} and reconstruction theory as in e.g. 
\cite[Theorem 5.4.1]{EtGeNiOs-tensor-categories} provide a 
cotriangular Hopf algebra $\Gamma_{E}=H(F_{E})$ for $E\in GL_{m}(\mC)$ 
as above such that, as monoidal categories,
\begin{gather*}
\coRep_{\BC}\Gamma_{E}\cong\cD
\end{gather*}
and such that $F_{E}$ becomes the forgetful functor.

It remains to argue that $\mathrm{tr}(E^{T}E^{-1})=-2$ admits a solution for every $m\in\BN_{\geq 2}$. Indeed, we can take
\begin{gather*}
E=
\left(\begin{array}{@{}c|c@{}}
id_{m-2} & 0
\\
\hline
0 & \begin{matrix}0 & x\\-1 & 0\end{matrix}
\end{array}\right)
,
\end{gather*}
for $x$ a solution to the equation $x^{2}-x(m-1)+2=0$ if $m\neq 2$ (which always has two solutions for $m\in\BN_{\geq 3}$), and $x=1$ for $m=2$. 
\end{proof}

\autoref{t:mainfalse} implies \autoref{t2:main}.(b) and 
\autoref{t3:main}.(c).

\begin{Remark}\label{r:mainfalse}
There is also a quantum version of \autoref{t:mainfalse} 
where one replaces the trace condition by 
$q\mathrm{tr}(E^{T}E^{-1})=-1-q^{2}\in\bk$.
\end{Remark}

%%%%%%%%%%%%%%%%%%%%%%%%%%%%%%%%%%%%%%%%%%%%%%

\section{Questions}\label{s:questions}

%%%%%%%%%%%%%%%%%%%%%%%%%%%%%%%%%%%%%%%%%%%%%%

We list a few open questions regarding $b_{n}^{\Gamma,V}$ and $\beta^{\Gamma,V}$.

\begin{Question}\label{q:bounds}
Let $\Gamma_{2}=SL_{2}$ and let $V_{2}$ denote its vector representation. It is easy to observe that the number of summands in $V_{2}^{\otimes n}$ in positive characteristic is bounded by the corresponding number in characteristic zero. In particular,
\autoref{e:clt} 
and \autoref{l:sl2} imply that
for arbitrary $p$ we have 
\begin{gather*}
A\cdot
2^{n}/n^{\alpha}\leq
b_{n}^{\Gamma_{2},V_{2}}
\leq B\cdot
2^{n}/n^{1/2}
\end{gather*}
for $A,B\in\BR_{>0}$ and $\alpha=1+(\log_{2}p)^{-1}$.

So we ask: For fixed $p$, is there some $\delta=\delta(p)\in\BR_{>0}$ for which
\begin{gather*}
b_{n}^{\Gamma_{2},V_{2}}\in\Theta(2^{n}n^{-\delta})?
\end{gather*}
\end{Question}
Note that in characteristic zero we have $\delta=1/2$ by \autoref{e:clt}.

\begin{Remark}\label{r:bounds}
The following observation was communicated to us by Pavel Etingof.
For $p=2$ the value $\delta\in\BR_{>0}$ in \autoref{q:bounds}
appears to be $\delta=\frac{1}{2}\log_{2}(8/3)$, which is approximately $0.708$. For example,
\begin{gather*}
\begin{tikzpicture}[anchorbase]
\node at (0,0) {\includegraphics[height=5cm]{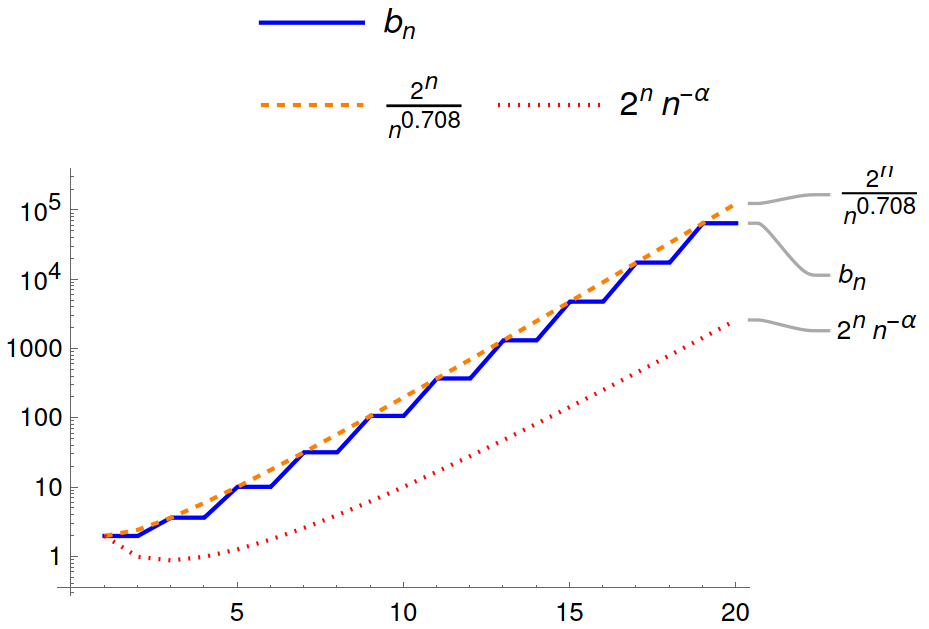}};
\end{tikzpicture}
,
\begin{tikzpicture}[anchorbase]
\node at (0,0) {\includegraphics[height=5cm]{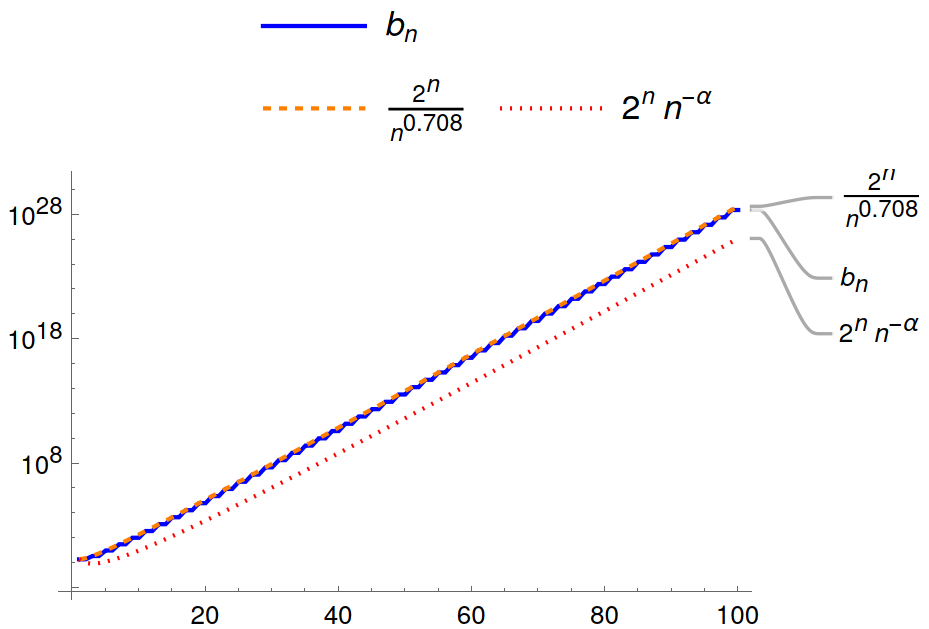}};
\end{tikzpicture}
\hspace*{-0.25cm}.
\end{gather*}
are Mathematica log plots for $p=2$ with $\alpha$ as in \autoref{l:sl2}.

The motivation for $\delta=\frac{1}{2}\log_{2}(8/3)$ is as follows. Consider the random variable $\big(\dim T(m)\big)^{s}$, where $T(m)$ as before denotes an indecomposable tilting $\Gamma_{2}$-representation with highest weight $m$ in $[n/2,n-1]$, and consider the uniform distribution of $m$. Then the expectation value 
$E\big(\big(\dim T(m)\big)^{s}\big)$ of this is, for $n\gg 0$, proportional to $n^{f(s)}$, where $f(s)=s-1+\log_{2}(1+2^{s})$. This can be proven using the character formula for $T(m)$ in \autoref{e:sl2}. In particular, the average of $1/\dim T(m)$ is $n^{f(-1)}=n^{-2+\log_{2}(3/2)}=n^{-\log_{2}(8/3)}$. Now, in characteristic zero the same type of argument would give that $n^{-1}$ is proportional to $E\big(1/\dim T(m)\big)$.
By \autoref{e:main} this suggests to take the square root, which gives the correct result up to a factor, i.e. in characteristic zero we have
\begin{gather*}
\sqrt{E\big(1/\dim T(m)\big)}=\sqrt{n^{-1}}=n^{-\delta}.
\end{gather*}
For $\mathrm{char}(\bk)=2$ this then suggests to take $\delta=\frac{1}{2}\log_{2}(8/3)\approx 0.708$ which, empirically speaking, seems to be correct, see the Mathematica output above.

A similar calculation could work for any prime. For example, for $p=3$ 
the above strategy gives $\delta=\frac{1}{2}\log_{2}(9/2)\approx 0.6845$
and
\begin{gather*}
\begin{tikzpicture}[anchorbase]
\node at (0,0) {\includegraphics[height=5cm]{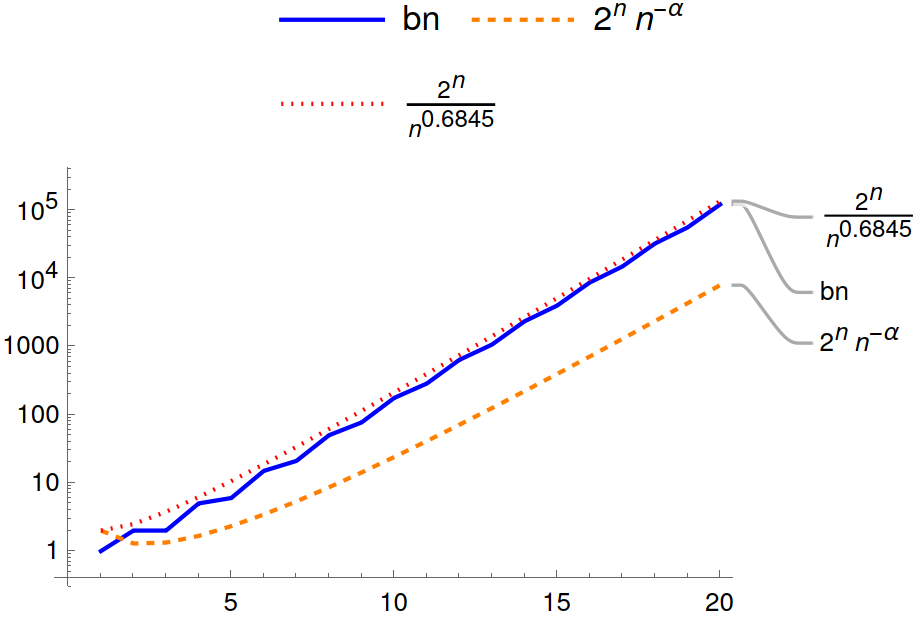}};
\end{tikzpicture}
,
\end{gather*}
is the associated Mathematica log plot.
\end{Remark}

More general than \autoref{q:bounds}, namely for all 
algebraic objects where a version of \autoref{t:main} holds, we could ask:

\begin{Question}\label{q:boundssecond}
For fixed $p$, is it true that
\begin{gather*}
b_{n}^{\Gamma,V}\in\Theta\big((\dim V)^{n}n^{-\delta}\big)
\end{gather*}
for some $\delta=\delta(p)\in\BR_{>0}$?
\end{Question}

For example even for $\Gamma_{3}=SL_{3}$ and $V_{3}$ its vector representation we do not know whether 
$b_{n}^{\Gamma_{3},V_{3}}$ is bounded from below by $A\cdot 3^{n}/n^{\delta}$ for some $A\in\BR_{>0}$.

\begin{Question}\label{q:boundsthird}
Assume that there exists $A\in\BR_{>0}$ such 
that 
\begin{gather*}
A\cdot(\dim V)^{n}\leq b_{n}^{\Gamma,V}?
\end{gather*}
What can we say about $\Gamma$ and $V$ where $\Gamma$ is assumed to be a group scheme?
\end{Question}

Recall that in the case $\mathrm{char}(\bk)=0$, \autoref{q:boundsthird} is answered in \autoref{p:fastgrowth}.
Note that the same answer cannot hold verbatim in positive characteristic (even after generalizing `torus' to group of multiplicative type; as we could also do in characteristic zero at no additional cost): there are finite group schemes which have a connected component which is not of multiplicative type, for instance infinitesimal group schemes.

\begin{Question}\label{q:tannakian}
Let $\Ver_{p}$ be the universal Verlinde category, see for instance \cite[Section 3]{Os-sym-fusion}, and $G$ an affine group scheme in $\Ver_{p}$ with a representation $X$ on an object $X_{0}\in\Ver_{p}$. Does \autoref{t2:main}(a) extend to the property that
\begin{gather*}
\beta^{G,X}=\mathrm{FPdim}(X_{0}),
\end{gather*}
with $\mathrm{FPdim}$ the Frobenius--Perron (sometimes called Perron--Frobenius) dimension. This formulation, with $\Ver_{p}$ replaced by its subcategories of (super) vector spaces corresponds precisely to \autoref{t2:main}(a) for (super)groups.
\end{Question}

Let $\cD$ be a Karoubian symmetric monoidal category with the property that the tensor product of two non-zero objects is not zero.
We say that
$\cD$ is of moderate decomposition growth
if for any object $X$ the sequence $b_{n}^{\Gamma,X}$ 
is bounded by a geometric progression. In this case the limit
$\beta^{\cD,X}$
exists, as before, by (a version of) Fekete's Subadditive lemma.

\begin{Question}\label{q:properties}
Does $\beta^{\cD,X}$ have nice properties? Is it additive? Is it multiplicative? Does it take values in algebraic numbers?
\end{Question} 

With respect to \autoref{RemSupP} we can ask:

\begin{Question}\label{q:super}
Does \autoref{t3:main}(c) extend to positive characteristic? In other words, is the number $b'_n\le b^{GL_{M|N}, V_{M|N}}_n$ of direct summands in $V_{M|N}^{\otimes n}$ which are forced to exist by considering idempotents in $\bk S_n$ (and the action of $\bk S_n$ on $V_{M|N}^{\otimes n}$) sufficient to ensure $\lim_{n\to\infty}\sqrt[n]{b'_n}=M+N$.
\end{Question}

Regarding the quantum case in \autoref{s:qgroups} one could ask:

\begin{Question}\label{q:quantum}
Let $\Gamma=U_{q}^{\bk}(\mathfrak{g})$ with $\mathfrak{g}$ not of type $A$, e.g. $\mathfrak{g}=\mathfrak{so}_{M}$. Do we still have $\beta^{\Gamma,V}=\dim V$?
\end{Question}

\section*{Declarations}

\subsection*{Ethical Approval} 
Not applicable.
 
\subsection*{Competing interests} 
No financial or personal competing interest.
 
\subsection*{Authors' contributions }
All authors contributed equally with respect to every section of the paper.
 
\subsection*{Funding} 
K.C. was partly supported by ARC grant DP200100712. D.T. was supported, in part, by the Australian Research Council.
 
\subsection*{Availability of data and materials} 
Not applicable.

\end{document}